\documentclass[a4paper]{article}

\usepackage[english]{babel}
\usepackage[utf8x]{inputenc}
\usepackage[T1]{fontenc}
\usepackage{lipsum}
\usepackage{blindtext}
\usepackage{graphicx,amssymb,amsmath,textcomp}
\usepackage{newpxtext} % ``New PX" font with original code \usepackage{newpxtext,newpxmath}
\usepackage{cancel}
\usepackage[a4paper,top=3cm,bottom=2cm,left=3cm,right=3cm,marginparwidth=1.75cm]{geometry}
\usepackage{tikz,tikz-cd}
\usetikzlibrary{matrix,arrows,positioning,calc}
\usepackage{verbatim}

\usepackage[nottoc]{tocbibind} %Includes "References" in the table of contents
\usepackage{amsmath,amsthm}
\usepackage{faktor}
\usepackage{xfrac} 
\usepackage{graphicx}
\usepackage[colorinlistoftodos]{todonotes}
\usepackage[colorlinks=true, allcolors=blue]{hyperref}
\usepackage{mathtools}

\newcommand{\Hom}{{\rm Hom}}
\newcommand{\vdim}{{\rm vdim \ }}

\newcommand{\rank}{{\rm rank \ }}
\newcommand{\cspec}{{\rm spec}}
\newcommand{\spec}{{\rm Spec}}

\newcommand{\N}{{\mathbb N}}

\newcommand{\R}{{\mathbb R}}

\newcommand{\K}{{\mathbb K}}
\newcommand{\tx}{\tilde{x}}

\newcommand{\dR}{{d_{dR}}}
\newcommand{\bfem}[1]{\textbf{\emph{#1}}}

\title{\textbf{Legendrian Structures in Derived Geometry}}
\author{Kadri İlker Berktav \footnote{Bilkent University, Department of Mathematics, Ankara, Turkey;  e-mail: kadri.berktav@bilkent.edu.tr. The author acknowledges that part of this research was carried out at the Institute of Mathematics, University of Zurich, under TUBITAK-2219 Grant Programme.}
} 

\date{\vspace{-5ex}}
\begin{document}
	
%%%%%%%%%%%%%%%%%%%% Text italic %%%%%%%%%%%%%%%%%%%%%%%%%%%%
\theoremstyle{plain}
\newtheorem{theorem}{Theorem}[section] %[section] or [subsection]
\newtheorem{lemma}[theorem]{Lemma} 
%NOT: araya [theorem] yazdigimiz zaman Lemma/propositon etc.'lari kendi icinde siralamak yerine 2. kod ile verilen genel numaralandirmayi (section veya subsection) takip eder.
\newtheorem{proposition}[theorem]{Proposition} %[section] or [subsection]
\newtheorem{corollary}[theorem]{Corollary}%[section]

%%%%%%%%%%%%%%%%%%%% Text roman %%%%%%%%%%%%%%%%%%%%%%%%%%%%%
\theoremstyle{definition}
\newtheorem{notations}[theorem]{Notations}
\newtheorem{notation}[theorem]{Notation}
\newtheorem{remark}[theorem]{Remark}
\newtheorem{observation}[theorem]{Observation}
\newtheorem{definition}[theorem]{Definition}
\newtheorem{condition}[theorem]{Condition}
\newtheorem{construction}[theorem]{Construction}
\newtheorem{example}[theorem]{Example}%[section]

\let\pf\proof
\let\epf\endproof
\numberwithin{equation}{section} %or {subsection}

\maketitle

\begin{abstract}
 This is the third installment in a series of papers \cite{kib,Berktav2024} on the subject of shifted contact structures on derived stacks.  In this paper,   we formally introduce the notion of a Legendrian structure in the derived context and provide natural constructions. We then present affine models and prove a Legendrian-Darboux theorem for the Legendrians in  contact derived schemes.
 \end{abstract}

\tableofcontents

    \section{Introduction}

    Derived geometry provides more generalized versions of well-known geometric structures, allowing for a study of their properties on higher spaces. For example, derived Symplectic and Poisson geometries have been introduced and examined in  \cite{PTVV,CPTVV}, along with several other applications and local constructions as in \cite{Brav,JS,BenBassat}.  
    
    In a series of papers \cite{kib,Berktav2024}, we have introduced \emph{derived contact geometry} and provided some results  inspired by the works \cite{Brav,BenBassat} and  the classical symplectic-contact dictionary. In brief, \cite{kib} defines \emph{shifted contact structures on derived stacks} and studies their local theory, providing a \emph{Darboux-type theorem}  and the notion of \emph{symplectification} for contact derived schemes; \cite{Berktav2024} extends the results of \cite{kib}  to derived Artin stacks and presents some further examples. Other attempts have also been made to develop the theory of shifted contact structures: we refer to \cite{Maglio2024}, investigating \emph{{$0$-shifted} and {$+1$-shifted contact structures} on differentiable stacks}.
   
  This paper, on the other hand,  examines shifted structures on the morphisms of  derived stacks.  In derived symplectic geometry, for example, \cite{PTVV} defines the notions of \emph{isotropic and Lagrangian structures} on a given morphism $f: {\bf Y}\rightarrow ({\bf X}, \omega_{\bf X})$ of derived stacks with shifted symplectic target. Regarding  local models for these structures, Joyce and Safronov \cite[Theorem 3.7]{JS} %showed that every Lagrangian $f: {\bf Y}\rightarrow ({\bf X}, \omega_{\bf X})$ in $k$-shifted symplectic derived scheme $ ({\bf X}, \omega_{\bf X}) $, with $k<0,$ is locally modeled on explicit ``Lagrangian-Darboux forms". This essentially 
    provide the \emph{derived version of the classical Lagrangian neighborhood theorem}, which inspires our current work.

 In the same spirit, this paper essentially aims to provide  the \emph{derived version of the classical Legendrian neighborhood theorem}.  Of course, one first requires analogous structures in derived contact geometry. To this end, the results and methods of \cite{JS} will be central for our constructions.
     
\paragraph{Results of the paper.}
In this sequel,  we formally introduce and study \emph{Legendrian structures} in the context of derived contact geometry. Let us overview our results.

%\begin{definition}(\emph{Informal})\label{defn_brief defn of LEG strc}
 	 In brief, we will define a \bfem{Legendrian structure} on a morphism $f: \bf Y \rightarrow X$ of derived Artin stacks, with a derived contact target $ (\mathcal{K}, \alpha, L[n])$, to be a \emph{non-degenerate isotropic structure} on the morphism $\mathbb{T}_{{\bf Y}} \rightarrow f^*(\mathcal{K})$. We then call $\bf Y$ (or $f$) a \emph{Legendrian} in $\bf X$. (See Definitions \ref{defn_isotropic for Leg} \& \ref{defn_LEG strc}.) %which will be formulated later. 
%\end{definition}

Once the proper definition and key observations are given, we prove the following results.
%We first give some interesting results including a natural example:
\begin{theorem}\label{thm_1point}
Any  Legendrian structure on the morphism ${\bf Y} \rightarrow \star_{n}$, where $\star_{n}$ denotes the point with its canonical $n$-shifted contact structure, induces an $(n-1)$-shifted exact symplectic structure on the derived scheme $\bf Y$ and an $(n-1)$-shifted contact structure on the product space ${\bf Y} \times \mathbb{A}^1[n-1]$, where $\mathbb{A}^1[n-1]$ is the \emph{$(n-1)$-shifted affine line} corepresented by the polynomial algebra on a variable in cohomological degree $ 1-n $. (cf.  Proposition \ref{prop_Leg on f with point target}.)
\end{theorem}

\begin{theorem}\label{thm_2}
	Let $ \mathrm{J^1[n]}{\bf X}$ be the \emph{$n$-shifted 1-jet stack} of $\bf X$ with the canonical $n$-shifted contact structure. Then the zero section $j: {\bf X} \rightarrow \mathrm{J^1[n]}{ \bf X}$ carries a natural Legendrian structure. (cf. Theorem \ref{thm_on zero section}.)
\end{theorem} 

%We then provide a Legendrian-Darboux-type theorem. In short, we have:

\begin{theorem} \label{THM3}
	 Legendrians in negatively shifted contact derived $\K$-schemes are locally modeled on standard Legendrian-Darboux forms (cf. Construction \ref{Leg model} and Theorem \ref{Legendrian-Darboux}).
\end{theorem}
    
    \paragraph{Organisation of the paper.}
    Section \ref{section_review of shifted symplectic strcs} outlines
    Lagrangian structures in derived symplectic geometry.
    In Section \ref{section_Leg strcs}, we formally define the notion of a \textit{Legendrian structure} on a morphism of derived stacks, with a shifted contact target, and present some natural results (the proofs of Theorems \ref{thm_1point} \& \ref{thm_2} above). Finally, Section \ref{section_local models and LEG darboux} provides a prototype construction for  such structures (cf. Construction \ref{Leg model}) and proves Theorem \ref{THM3} above (cf. Theorem \ref{Legendrian-Darboux}).
    
     \paragraph {Acknowledgments.} I thank Alberto Cattaneo and  Ödül Tetik for helpful conversations. I am grateful to the Institute of Mathematics, University of Zurich, where an early draft of this paper was prepared. %I personally benefited a lot from hospitality and research environment of the Institute. 
     The author acknowledges support of the Scientific and Technological Research Council of Turkey (T\"{U}BİTAK) under 2219-International Postdoctoral Research Fellowship (2021-1).
     
      Most parts of this paper were revisited and improved at the time when I joined Bilkent University Department of Mathematics in Fall 2023. Since then, I have benefited a lot from the discussions with Fırat Arıkan and Özgür Kişisel (both from Dept. of Math, Middle East Technical University). I am grateful to them for our research meetings and helpful suggestions.  
      
      I also wish to warmly thank the anonymous referee(s) for their valuable comments and suggestions, which helped a lot and improved the quality of the prequels \cite{kib,Berktav2024}, leading to several corrections and modifications for the current preprint. 

     \paragraph{Conventions.} Throughout the paper, $ \mathbb{K} $ will be an algebraically closed field of characteristic zero. All cdgas will be graded in nonpositive degrees and  over $ \mathbb{K}.$ All classical $ \K $-schemes are assumed \emph{locally of finite type}, and all derived $ \K $-schemes/stacks $ \textbf{X} $ are assumed to be \emph{locally finitely presented,} by which we mean that there exits a cover of $\bf X$ by affine opens $\spec  A$ with $A$ a finitely generated graded algebra. 
     
    \section{Recollection: Lagrangians in the derived context} \label{section_review of shifted symplectic strcs}

\subsection{Derived symplectic geometry}   \label{secion_PTVV's symplectic geometry on spec A} 
  Pantev et al. \cite{PTVV} define the \emph{simplicial sets of $p$-forms of degree $k$  and closed $p$-forms of degree $k$}  on derived stacks. Denote these simplicial sets by $\mathcal{A}^p({\bf X},k)$ and $\mathcal{A}^{p,cl}({\bf X},k)$, respectively.
  
   These definitions are in fact given first for affine derived $\K$-schemes. Later, both concepts are defined for  general derived stacks $\bf X$ in terms of mapping stacks $\mathcal{A}^p(-,k)$ and $\mathcal{A}^{p,cl}(-,k)$, respectively. % A summary of key ideas can be found in \cite[Section 3.4]{Brav}
  For a summary of some key ideas, see Appendix \ref{section_prelim}.

  Recall that any derived stack $\textbf{X}$ has a graded mixed \bfem{de Rham complex}, written $DR(\textbf{X})$, with $d_{tot}=d+ \dR$.
  Then given any derived stack $\textbf{X}$, we define \begin{equation}\label{defn_DR(X) as inf limit}
  	DR(\textbf{X})= \lim_{A \in cdga_{\K}^{\leq 0}, \ \spec A \rightarrow X} DR(A).
  \end{equation}  If, in addition, $\bf X$ is Artin and $\mathbb{L}_{\bf X}$ denotes its \emph{cotangent complex}, then we also have 
  \begin{equation}
  	DR(\textbf{X}) \simeq \Gamma(\textbf{X}, \mathrm{Sym}_{\mathcal{O}_{\textbf{X}}} (\mathbb{L}_{\textbf{X}}[1])).
  \end{equation}

  Let $ \textbf{X} = \spec A$ with $A$ a standard form cdga\footnote{It should be noted that the results that are cited or to be proven in this section are all about the \emph{local structure} of derived schemes. Thus, it is enough to consider the (refined) affine case.}, then we can take\footnote{Not true for an arbitrary cdga.} $ \Lambda^p\mathbb{L}_A=\Lambda^p \Omega^1_{A}$,  where $\Omega^1_{A}$ denotes the  \bfem{$A$-module  of K\"{a}hler differentials}.
  Then  we obtain \begin{equation}
  	DR(\textbf{X})= DR(A)\simeq\mathrm{Sym}_A (\Omega^1_A [1])\simeq\bigoplus \limits_{p\geq0}  \Lambda^p \Omega^1_{A} [p],
  \end{equation} where the graded mixed differential is given by the universal derivation $A\rightarrow \Omega^1_A$ extended via the Leibniz rule. Here, any element in $  (\Lambda^p \Omega^1_{A})^k [p] $ has \bfem{(cohomological) degree} $k-p$ and \bfem{weight} $p.$
  As our local models can always be formed by such cdgas \cite{BenBassat}, we will simply  work with  the  module $\Omega^1_{A}$ instead of $ \mathbb{L}_A$.  
  
  Therefore,  elements of $\mathcal{A}^p({\bf X},k)$  form a simplicial set such that $k$-cohomology classes of the complex $ \big(\Lambda^p \Omega^1_{A}, d\big) $  correspond to the connected components of this simplicial set. Likewise,  the connected components of $\mathcal{A}^{p,cl}({\bf X},k)$ are identified with the $k$-cohomology classes of the complex $\prod_{i\geq0} \big(\Lambda^{p+i} \Omega^1_{A}[-i], d_{tot}\big)$.   Then we have the following definitions.%We want to work with explicit representatives for these cohomology classes. 

  \begin{definition} \label{defn_p form of deg k}
  	Let $\textbf{X}=\spec A$ be an affine derived $\mathbb{K}$-scheme for $A$ a minimal standard form cdga. A \emph{\textbf{$k$-shifted $p$-form on $\bf X$}} for $p\geq 0$ and $k \leq 0$ is an element 
  	$ \omega^0 \in \big(\Lambda^p \Omega^1_{A}\big)^k \text{ with } d\omega^0=0.  $
  	
  	Two $ p $-forms  $\omega_1^0, \omega_2^0 $ of degrees $k$ are \bfem{equivalent}, written $\omega_1^0 \sim  \omega_2^0$, if $\exists \alpha_{1,2} \in \big(\Lambda^p \Omega^1_{A}\big)^{k-1}$, a \bfem{path}, s.t. $\omega_1^0-\omega_2^0= d\alpha_{1,2}$.
  \end{definition} %Note that an element  $\omega^0$ defines a cohomology class $[\omega^0] \in H^k\big(\Lambda^p \Omega^1_{A}, d\big),$ where 

  \begin{definition}\label{defn_closed p form}
  	Let $\textbf{X}=\spec A$ be an affine derived $\mathbb{K}$-scheme with $A$ a minimal standard form cdga. A \emph{\textbf{closed $k$-shifted $p$-form on $\bf X$}} for $p\geq 0$ and $k \leq 0$ is a sequence $\omega=(\omega^0, \omega^1, \dots)$ with $\omega^i \in \big(\Lambda^{p+i} \Omega^1_{A}\big)^{k-i}$ such that $d_{tot}\omega=0$, which splits according to weights as $d\omega^0=0$ in $ \big(\Lambda^p \Omega^1_{A}\big)^{k+1}$,
  	and $d_{dR}\omega^i + d\omega^{i+1}=0$ in $\big(\Lambda^{p+i+1} \Omega^1_{A}\big)^{k-i}$, $i \geq 0.$ 
  	
  		We say two closed $ p $-forms  $\omega=(\omega^0, \omega^1, \dots), \sigma=(\sigma^0, \sigma^1, \dots) $ of degrees $k$ are \bfem{equivalent}, written $\omega \sim \sigma$, if there exists a sequence (a \bfem{path}) $\alpha=(\alpha^0,\alpha^1,\dots)$ with $\alpha^i \in (\Lambda^{p+i} \Omega_A^1)^{k-i-1}$ for $i=0,1,\dots,$ satisfying \begin{equation*}
  		\omega^0-\sigma^0=d\alpha^0 \text{ and } \omega^{i+1} + \sigma^{i+1}= d_{dR}\alpha^i + d\alpha^{i+1}.
  	\end{equation*}
  	
\end{definition} Observe that a closed $k$-shifted $p$-form consists of an actual $k$-shifted $p$-form $ \omega^0 $ and the data $ (\omega^i)_{i>0} $ of $ \omega^0 $ being coherently $d_{dR}$-closed up to homotopy.
It then follows that there also exists a natural projection morphism  
$ \pi: \mathcal{A}^{p,cl}(\textbf{X},k)\longrightarrow \mathcal{A}^{p}(\textbf{X},k), \ \ \omega=(\omega^i)_{i\geq 0} \longmapsto \omega^0. $
  
 Note also that the simplicial set $\mathcal{A}^{p,cl}({\textbf{X}},k)$ is usually more complicated than $ \mathcal{A}^p({\textbf{X}},k)$ even if we make nice assumptions on a general derive stack $\textbf{X}$. Thanks to \cite[Prop.  1.14]{PTVV}, we have the following identification for the space $ \mathcal{A}^p({X},k)$, making it more tractable: \begin{equation}
 	\mathcal{A}^p(\textbf{X},k)\simeq Map_{QCoh(\textbf{X})}(\mathcal{O}_{\textbf{X}}, \wedge^p\mathbb{L}_{\textbf{X}}[k]). 
 \end{equation}
 Thus, any 2-form of degree $k$  induces a morphism $\mathcal{O}_{\textbf{X}} \rightarrow  \wedge^2\mathbb{L}_{\textbf{X}}[k]$ in $QCoh(\textbf{X})$, and hence, by duality, a morphism $\mathbb{T}_{\textbf{X}} \wedge \mathbb{T}_{\textbf{X}} \rightarrow \mathcal{O}_{\textbf{X}}[k]$. By adjunction, this gives the induced morphism $ \mathbb{T}_{\textbf{X}} \rightarrow \mathbb{L}_{\textbf{X}} [k]$, leading to the following definition:
  
  \begin{definition}
  	A closed $k$-shifted $2$-form $ \omega=(\omega^i)_{i\geq 0} $ on  $\textbf{X}=\spec A$ for a (minimal) standard form cdga $A$ is called a \emph{\textbf{$k$-shifted symplectic structure}} if  the induced map $$ \omega^0\cdot: \mathbb{T}_{A} \rightarrow \Omega^1_{A}[k], \ Y \mapsto \iota_{Y} \omega^0,$$ is a quasi-isomorphism, where $\mathbb{T}_{A}=(\mathbb{L}_{A})^{\vee}\simeq\Hom_{A}(\Omega^1_{A},A)$\footnote{Thanks to the identification $\mathbb{L}_A\simeq (\Omega_A^1, d)$ for $A$ a (minimal) standard form cdga.} is the \textit{tangent complex of $A$.}
  \end{definition}

\paragraph{Shifted symplectic Darboux models for derived schemes.} \label{the pair}
One of the main theorems in  \cite{Brav} provides the shifted version of the classical Darboux theorem in symplectic geometry. The statement is as follows. 

\begin{theorem} (\cite[Theorem 5.18]{Brav}) \label{Symplectic darboux}
	Given a derived $\mathbb{K}$-scheme $X$ with a $k$-shifted symplectic form $\omega'$ for $k<0$ and $x\in X$, there is a local model $\big(A, f: \spec A \hookrightarrow X, \omega \big)$  and $p \in \spec (H^0(A))$ such that $f$ is an open inclusion with $f(p)=x$, $A$ is a standard form that is minimal at  $p$, and $\omega$ is a $k$-shifted symplectic form on $\spec A$ such that $A, \omega$ are in Darboux form, and $f^*(\omega') \sim \omega$  in the space of $k$-shifted closed 2-forms.
\end{theorem} 

In fact, it has been proven in \cite[Theorem 5.18]{Brav} that such $ \omega $  can be constructed explicitly depending on the integer $k<0$. Indeed, there are three cases in total: 
 \emph{ $(1) \ k  $  is odd;  $ (2) \ k/2 $  is even;  and  $ (3) \ k/2  $ is odd.} % Therefore, Theorem \ref{Symplectic darboux} means that every $k$-shifted symlectic derived $\K$-scheme $(\bf X, \omega')$ is Zariski locally equivalent to $(\spec A, \omega)$ for some $A,\omega$ in Darboux form, given by one of the  cases above. 
 For instance, when $k$ is odd, one can find a minimal standard form cdga $A$, with $ `` $coordinates" $x^{-i}_j, y^{k+i}_j \in A$, and a Zariski open inclusion $f: \spec A \hookrightarrow \bf X$ so that $f^*(\omega') \sim \omega=(\omega^0, 0, 0, \dots)$ and
$  \omega^0= \sum_{i,j} d_{dR}x_j^{-i} d_{dR}y_j^{k+i}.$ We will not give any further detail on the aforementioned cases in this paper. Instead, we refer to \cite[Examples 5.8, 5.9, and 5.10]{Brav}.

\subsection{Lagrangians in symplectic derived stacks} \label{section_Lag in PTVV}
 
In classical  symplectic  geometry, there are some natural subobjects of symplectic manifolds, called \emph{Lagrangians}\footnote{A \bfem{Lagrangian} $\iota: L\hookrightarrow (M,\omega)$ is a submanifold of dimension $\frac{1}{2}\dim M$ such that $\iota^*\omega\equiv0$.}. There is a well-known result on the standard model for tubular neighborhoods of Lagrangian submanifolds: \emph{Weinstein's Lagrangian neighborhood theorem} states \cite{McDuff2017} that given a Lagrangian submanifold $L\subset (M, \omega)$, a neighborhood of the zero section of its cotangent bundle $T^*L$ provides a universal model for the neighborhood of $L$ itself. That is, we have: \begin{theorem}
	Every Lagrangian $L\subset (M, \omega)$ has a tubular neighborhood symplectomorphic to a neighborhood of the zero section of $T^*L$.
\end{theorem}%In what follows, we shall discuss  analogous results in derived symplectic geometry.

Joyce and Safronov \cite{JS} provide an analogous result in the context of derived symplectic geometry. It should be noted that  PTVV's paper \cite{PTVV} introduces the notions of \emph{isotropic} and \emph{Lagrangian structures} on a given morphism $f: \mathcal{L}\rightarrow ({\bf X}, \omega_{\bf X})$ in $dStk$  to study, for instance, derived symplectic structures on fiber products of derived stacks and that on mapping stacks \cite[Theorems 2.9 \& 2.5]{PTVV}. Now, let us discuss the notions of interest in detail. We follow \cite{PTVV, JS}.
%In this regard, we have the following important result:
%\begin{theorem} \cite[Theorem 2.9.]{PTVV}
%	Let $L_1\xrightarrow{f}({\bf X}, \omega_{\bf X}) \xleftarrow{g} L_2$ be a diagram of derived Artin $\K$-stacks, where $\omega_{\bf X}$ is a $k$-shifted symplectic structure and $f,g$ are morphisms with Lagrangian structures, then the fiber product $L_1 \times_{\bf X} L_2$ admits a canonical $(k-1)$-shifted symplectic structure. 
%	\end{theorem}

%Now, let us discuss the notions above in detail. We follow \cite{PTVV, JS}.
Denote by $\mathcal{A}^p(-,n), \mathcal{A}^{p,cl}(-,n)$ the \emph{derived stacks of $ p $-forms of degree $ n $} and \emph{closed $ p $-forms of degree $ n $}
as introduced in the previous section. %Then by construction, we also have equivalences \begin{equation} \label{equivalences_loop space vs form}
%	\mathcal{A}^p(n)\simeq \Omega \mathcal{A}^p(n+1) \text{ and }  	\mathcal{A}^{p,cl}(n)\simeq \Omega \mathcal{A}^{p,cl}(n+1).
%\end{equation}
\begin{definition} \label{defn_isotropic Lag}
Let $ ({\bf X}, \omega_{\bf X}) $ be a $k$-shifted symplectic derived Artin stack and $f: \mathcal{L}\rightarrow ({\bf X}, \omega_{\bf X})$ a morphism of derived Artin stacks.
By an \emph{\textbf{isotropic structure}} on $f$ (relative to $ \omega_{\bf X} $), we mean a homotopy $h_{\mathcal{L}}$ from $0$ to $f^*(\omega_{\bf X})$ in the simplicial set $\mathcal{A}^{2, cl} (\mathcal{L}, k)$; i.e., a nullhomotopy of $f^*(\omega_{\bf X})$.  Denote the space of such structures by 
\begin{equation*}
Isot(f, \omega_{\bf X}) := Path_{\mathcal{A}^{2, cl}_{\K} (\mathcal{L}, k)} {\big(0, f^*(\omega_{\bf X})\big) } .
\end{equation*}
\end{definition}

%\begin{remark} \label{rmk_nulhomotopy as a form}
%	Note that if a nullhomotopy $h'_{\mathcal{L}}$ of an element in $\mathcal{A}^{2 (,cl)}(k)$ defines a loop in the space of (closed) 2-forms of degree $ k $, from the equivalences (\ref{equivalences_loop space vs form}), one can equivalently see $h'_{\mathcal{L}}$ as a (closed) 2-form of degree $(k-1).$ 
%\end{remark}

\begin{observation}

Let $ h_{\mathcal{L}}\in Isot(f, \omega_{\bf X}) $, then using the natural projection $$ \pi: \mathcal{A}^{p,cl}({\bf X},k)\longrightarrow \mathcal{A}^{p}({\bf X},k), \ \ \omega_{\bf X}=(\omega^i)_{i\geq 0} \longmapsto \omega^0,$$ we can obtain a homotopy $ h_{\mathcal{L}}^0 $ from $0$ to $f^*(\omega^0_{\bf X})$ in $ \mathcal{A}^{2} (\mathcal{L}, k),$ which induces the following 2-commutative diagram (via $h_{\mathcal{L}}^0 \cdot$), along with the induced morphism $\chi_{ h_{\mathcal{L}}}: 	\mathbb{T}_{\mathcal{L}/\bf X} \rightarrow \mathbb{L}_{\mathcal{L}}[k-1]:$

\begin{equation} \label{isotropic strc}
\begin{tikzpicture}
\matrix (m) [matrix of math nodes,row sep=1.5em,column sep=4.5em,minimum width=1.5 em] {
	\mathbb{T}_{\mathcal{L}}   & { } & 0 \\
f^*(\mathbb{T}_{\bf{X}})	& f^*(\mathbb{L}_{\bf {X}})[k] &  \mathbb{L}_{\mathcal{L}}[k], \\
	\mathbb{T}_{\mathcal{L}/\bf X}[1] &  \\};
\path[-stealth]
(m-1-1) edge  node [left] { {\small $ \mathbb{L}_f^* $}} (m-2-1)
edge  node [above] { } (m-1-3)
(m-1-3) edge  node [right] { } (m-2-3)
(m-2-1) edge  node [below] {} node [above] {{\small $ f^*(\omega^0_{\bf X})\cdot  $}} (m-2-2)
%(m-1-1) edge  node [below] {} node [below] {{\small  higher stacks}} (m-3-2)
(m-2-2) edge  node  [above] {{\small $ \mathbb{L}_f [k] $}} (m-2-3)
(m-3-1) edge [dashed]  node [below] { $ \chi_{ h_{\mathcal{L}}} [1] $ } (m-2-3)
(m-2-1) edge  node [right] { } (m-3-1)
(m-1-2) edge [double]  node [left] {$ h_{\mathcal{L}}^0 \cdot $ } (m-2-2);
%edge [dashed,-] (m-2-1);
\end{tikzpicture}
\end{equation}
Here $ \mathbb{L}_f [k] $ is the shifted version of the moprhism in the  triangle $f^*(\mathbb{L}_{{\bf X}}) \rightarrow \mathbb{L}_{\mathcal{L}} \rightarrow  \mathbb{L}_{\mathcal{L}/\bf X}$ for  $f: \mathcal{L}\rightarrow ({\bf X}, \omega_{\bf X})$. The morphism $ f^*(\omega^0_{\bf X})\cdot  $ is the pullback of the induced map $ \omega_{\bf X}^0 \cdot :  \mathbb{T}_{\bf{X}} \xrightarrow {\sim} \mathbb{L}_{\bf {X}}[k]$. The vertical maps on the left hand side of the diagram are just the duals of ones in the canonical distinguished triangle above. The map $ \chi_{ h_{\mathcal{L}}} $ is described below, see Observation \ref{observation_gettting LAG non-deg }. 

\begin{observation} \label{observation_gettting LAG non-deg }
	The homotopy $ h_{\mathcal{L}}^0$ between 0 and $f^*(\omega_{\bf X}^0): f^*(\mathbb{T}_{\bf{X}}) \wedge f^*(\mathbb{T}_{\bf{X}}) \rightarrow \mathcal{O}_{\mathcal{L}} [k]$ induces a homotopy $ h' $ from 0 to the map $\mathbb{T}_{\mathcal{L}/\bf X} \otimes \mathbb{T}_{\mathcal{L}} \rightarrow \mathcal{O}_{\mathcal{L}} [k]$ obtained by the composition\footnote{The first two maps come from the exact triangle $\mathbb{T}_{\mathcal{L}/\bf X}\rightarrow \mathbb{T}_{\mathcal{L}} \rightarrow f^*(\mathbb{T}_{\bf{X}}).$} \begin{equation}
	\mathbb{T}_{\mathcal{L}/\bf X} \otimes \mathbb{T}_{\mathcal{L}} \rightarrow \mathbb{T}_{\mathcal{L}}\wedge \mathbb{T}_{\mathcal{L}} \rightarrow f^*(\mathbb{T}_{\bf{X}}) \wedge f^*(\mathbb{T}_{\bf{X}}) \rightarrow \mathcal{O}_{\mathcal{L}} [k].
\end{equation} On the other hand,  $ \mathbb{T}_{\mathcal{L}/\bf X} \rightarrow f^*(\mathbb{T}_{\bf{X}}) $ comes with a canonical homotopy\footnote{Thanks to the exact triangle $\mathbb{T}_{\mathcal{L}/\bf X}\rightarrow \mathbb{T}_{\mathcal{L}} \rightarrow f^*(\mathbb{T}_{\bf{X}}) $} to 0, and hence we obtain another induced homotopy $h''$ between 0 and the composition $\mathbb{T}_{\mathcal{L}/\bf X} \otimes \mathbb{T}_{\mathcal{L}} \rightarrow \mathcal{O}_{\mathcal{L}} [k]$. 

Combining the homotopies $h', h''$, we get a loop at 0 in $Map_{QCoh({\bf X})} \big(\mathbb{T}_{\mathcal{L}/\bf X} \otimes \mathbb{T}_{\mathcal{L}}, \mathcal{O}_{\mathcal{L}} [k]\big),$ which defines an element in  $$\pi_1\big(Map_{QCoh({\bf X})} (\mathbb{T}_{\mathcal{L}/\bf X} \otimes \mathbb{T}_{\mathcal{L}}, \mathcal{O}_{\mathcal{L}} [k])\big) \simeq \pi_0\big(Map_{QCoh({\bf X})} (\mathbb{T}_{\mathcal{L}/\bf X} \otimes \mathbb{T}_{\mathcal{L}}, \mathcal{O}_{\mathcal{L}} [k-1])\big).$$  Using adjunction, we get a morphism $\chi_{h_{\mathcal{L}}}: 	\mathbb{T}_{\mathcal{L}/\bf X} \rightarrow \mathbb{L}_{\mathcal{L}}[k-1]$ of perfect complexes. For more details, see \cite[\S 2.2]{PTVV}. 
\end{observation}

These observations lead to the definition of a \emph{Lagrangian structure} on the morphism $ f $ (with respect to $ \omega_{\bf X} $), which will simply be an \emph{isotropic
structure satisfying a non-degeneracy condition.} More precisely, we have:

\end{observation}
\begin{definition}
%	We say that $ h_{\mathcal{L}}^0 $ is \emph{non-degenerate} if the morphism $\chi_{h_{\mathcal{L}}}: 	\mathbb{T}_{\mathcal{L}/\bf X} \rightarrow \mathbb{L}_{\mathcal{L}}[k-1]$ in (\ref{isotropic strc}) is a quasi-isomorphism. 
	
	Let $f: \mathcal{L}\rightarrow ({\bf X}, \omega_{\bf X})$ be a morphism of derived Artin stacks with $k$-shifted symplectic target. An isotropic structure $ h_{\mathcal{L}} $ on  $f$ is   called \bfem{Lagrangian} if the induced morphism $\chi_{h_{\mathcal{L}}}: 	\mathbb{T}_{\mathcal{L}/\bf X} \rightarrow \mathbb{L}_{\mathcal{L}}[k-1]$ in (\ref{isotropic strc}) is a quasi-isomorphism (the \emph{non-degeneracy condition})\footnote{Equivalently, we can require the sequence $\mathbb{T}_{\mathcal{L}}\rightarrow f^*(\mathbb{T}_{{\bf X}}) \rightarrow \mathbb{L}_{\mathcal{L}}[k]$ to be a homotopy fiber sequence.}. In that case, we  say that $\mathcal{L}$ is \emph{Lagrangian} in $ ({\bf X}, \omega_{\bf X}). $
\end{definition}

\paragraph{Local description of isotropic structures.}
In the presence of ``nice" local models representing a given morphism $f: \mathcal{L}\rightarrow ({\bf X}, \omega_{\bf X})$, one can express the isotropy condition on $f$ in terms of certain local equations. For the meaning of ``nice" local models, see Appendix \ref{section_prelim}.

Let ${\bf X} \simeq \spec A$ and  $\mathcal{L} \simeq \spec B$ be affine derived $\K$-schemes, with $A,B$ some (standard form)\footnote{Thus, we can identify their cotangent complexes with the corresponding modules of K\"{a}hler differentials.} cdgas, and $f: \mathcal{L}\rightarrow ({\bf X}, \omega_{\bf X})$ be a morphism in $\mathrm{dAff}_{\K}$ with $k$-shifted symplectic target such that $f$ is induced by a morphism\footnote{That is, $f\simeq \spec \tilde{f}$} of cdgas $\tilde{f}:A \rightarrow B$. Then $\omega_{\bf X}$ lifts to a sequence $\omega_A=(\omega^0, \omega^1, \dots) \in \big(\prod_{i\geq0} \Lambda^{2+i} \Omega_A^1[-i], d_{tot}\big)$ such that $ d_{tot}\omega_A=0. $ 

Let $h: 0 \rightsquigarrow f^*\omega_A$ be an isotropic structure   on $f$,  then $[f^*\omega_A]=[0]$ in $H^k\big(\prod_{i\geq0} \Lambda^{2+i} \Omega_B^1[-i],  d_{tot}\big)$, and hence,  by definition, there is a sequence $ (h^0, h^1, \dots) $
 with $h^i\in \big(\Lambda^{2+i} \Omega^1_B\big)^{k-1-i}$ for $i=0,1,\dots,$  satisfying %that the isotropy condition can be expressed via the equations (see also \cite{JS})
 \begin{equation}
 dh^0= \tilde{f}_* (\omega^0), \ \ \ \ \ d_{dR}h^i + dh^{i+1}= \tilde{f}_*(\omega^{i+1}), \ \ i=0,1,\dots,
 \end{equation} where $ \tilde{f}_*: \Omega^1_A \rightarrow \Omega^1_B $ is the morphism induced by $\tilde{f}.$ Thus, the isotropic structure $h$  can equivalently be viewed as the sequence $ (h^0, h^1, \dots) $ above.

Recall from \cite[Prop. 5.7]{Brav} that since our local model is nice enough (given by standard form cdgas), we may take $\omega^0$ to be exact and $\omega^i=0$ for $i>0$, and write $\omega_A=(\omega^0, 0, \dots)$. In this case, the defining equations above reduce to
 \begin{equation}
 dh^0= \tilde{f}_* (\omega^0), \ \ \ \ \ d_{dR}h^i + dh^{i+1}= 0, \ \ i=0,1,\dots,
 \end{equation} 
 \begin{remark} \label{simplification for h}
 	As in the case of shifted symplectic structures with nice local models, one can also simplify the form of the isotropic structure $ (h^0, h^1, h^2, \dots) $ by using certain vanishing results from cyclic homology theory. In that respect,  \cite[Proposition 4.1]{JS} shows that, up to equivalence, $h^0$ can be taken to be exact and $h^i=0$ for $i>0$. Therefore, for an isotropic structure of the form $ (h^0, 0, 0, \dots) $, the defining equations  even further reduce to
 \begin{equation}
 dh^0= \tilde{f}_* (\omega^0), \ \ \ \ \ d_{dR}h^0=0.
 \end{equation} 
 \end{remark}
 \begin{observation}
 	When $({\bf X}, \omega_{\bf X})=(*, 0)$, i.e. $ \omega^0=0 $ as well, we obtain $ h:=(h^0, 0, 0, \dots) $, with $ dh^0= 0, \  d_{dR}h^0=0, $ and hence an element in $ \mathcal{A}^{2,cl}(\spec B, k-1).$ That is, we have a $(k-1)$-shifted pre-symplectic structure $h$ on $ \spec B.$ In addition, if the isotropic structure is non-degenerate, i.e. $\chi_h:\mathbb{T}_f \rightarrow \mathbb{L}_B[k-1]$ is a quasi-isomorphism, then the map $\chi_h$ reduces to  the non-degenerate contraction map $\chi_h\cdot:\mathbb{T}_B \rightarrow \mathbb{L}_B[k-1]$. Thus, $h$ defines a $(k-1)$-shifted symplectic structure on $ \spec B.$ In general, we have the following result.
 	\begin{proposition}
 		Let $ ({\bf X}, \omega_{\bf X})=(*, 0) $ with a trivial $k$-shifted symplectic structure $\omega_{\bf X}=0.$ Lagrangians $Y$ in $ (*, 0) $ are equivalent to $(k-1)$-shifted symplectic derived schemes $(Y, \omega_Y)$, where $\omega_Y$ is defined by the non-degenerate $k$-shifted isotropic structure on $Y\rightarrow (*, 0).$
 	\end{proposition}
 \end{observation}

\paragraph{A Lagrangian neighborhood theorem.} Before stating a Lagrangian Darboux-type theorem, we wish to present a simplification result for an isotropic structure $h=(h^1, h^2, \dots)$ on  $\spec \beta: \spec B \rightarrow (\spec A, \omega_{can})$ which was  briefly mentioned in Remark \ref{simplification for h}. The following result is central to prove the shifted version of the Lagrangian-Darboux theorem in the classical setup.

\begin{lemma} \cite[Proposition 4.1.]{JS} \label{lemma_Lag simplification}
	Let $A, \omega_{can}, B, \beta$ be as above, and $ h=(h^0, h^1, \dots) $
	 an isotropic structure for $\spec \beta$ satisfying $ dh^0= \beta_* (\omega^0), \  d_{dR}h^i + dh^{i+1}= 0, \ \ i\geq0. $	 
	 Then there exist elements $G'\in B^k$ and $\psi\in (\Omega_B^1)^{k-1}$ satisfying the equations $$ dG'=-\beta(H+H_+) \text{ and } d_{dR}G'+d\psi=-\beta_*(\phi+\phi_+) $$ such that the isotropic structure $ h=(h^0, h^1, \dots) $ is homotopic to $ \frac{1}{k-1} (d_{dR}\psi, 0, 0, \dots). $
\end{lemma}

Now, we state the ``derived" Lagrangian-Darboux theorem. For details on standard local models, see Appendix \ref{section_models for LAGNBH thm}.
\begin{theorem} \cite[Theorem 3.7.]{JS}
	Lagrangians $f:\mathcal{L} \rightarrow ({\bf X}, \omega_{\bf X})$ in $k$-shifted symplectic derived schemes, with $k<0$, are locally modeled on standard Lagrangian-Darboux forms as in Definition \ref{defn_Lag Darboux forms} and Construction \ref{Lag model}, with the following homotopy commutative diagram 
	\begin{equation} 
		\begin{tikzcd}
		\spec B \arrow[d, "\spec \beta"'] \arrow[r, "j", hook] & \mathcal{L} \arrow[d, "f"] \\
		\spec A \arrow[r, "\iota"', hook]                          & \bf X               
	\end{tikzcd}
	\end{equation} and the pullback $j^*(h_{\mathcal{L}})$ of the Lagrangian structure $h_{\mathcal{L}}$ on $f$ to $ \spec \beta $ so that $j^*(h_{\mathcal{L}})\simeq (h^0, 0, 0, \dots)$.% and $  h^0= \sum_{i,j} d_{dR}u_j^{-i} d_{dR}v_j^{k-1+i}.$
\end{theorem}

\section{Legendrian structures in the derived context} \label{section_Leg strcs}

\subsection{Derived contact structures} \label{section_shifted contact structures and Darboux forms}

In classical contact geometry, by a \emph{contact structure} on a smooth manifold $M^{2n+1}$, we mean a smooth field of tangent hyperplanes $ \xi \subset TM$ (of $ \rank 2n $) with the property that for any smooth locally defining 1–form $\alpha$, i.e. $\xi= \ker (\alpha)$, the 2-form $d_{dR}\alpha |_{\xi}$ is non-degenerate. Recall that contact manifolds are viewed as the odd-dimensional analogues of symplectic manifolds. In that respect, both structures have several common features; some of which are: a Darboux theorem, no local invariants, and  having interesting subobjects. For more details, we refer to \cite{Geiges}. 

In derived geometry, on the other hand,  \cite{kib,Berktav2024} introduce \textit{shifted contact structures} on  derived  stacks  and study their properties. %, with $n\in \Z,$ as a structure that consist of a  morphism of $\mathcal{K}\rightarrow \mathbb{T}_{{\bf X}}$ of perfect complexes, a line bundle $L$, and a locally defined $n$-shifted 1-form $\alpha: \mathbb{T}_{{\bf X}} \rightarrow \mathcal{O}_{{\bf X}}[k]$ satisfying a non-degeneracy condition. Moreover, \cite{kib,Berktav2024} present derived versions of some of the classical results in contact geometry.
Let us start  with some terminology and  results from \cite{kib,Berktav2024}.
\begin{definition}\label{defn_shiftedcontact} Let $\bf X$ be a locally finitely presented derived (Artin) stack.
	An \emph{\textbf{$ n $-shifted contact structure}} on $\bf X$ consists of \begin{itemize}
			\item a line bundle $L$, a perfect complex $\mathcal{K}$, and a morphism  $ \alpha: \mathbb{T}_{{\bf X}}\rightarrow L[n] $ %and a perfect complex $\mathcal{K}$ on $\bf X$ 
			fitting into the fiber-cofiber sequence  $$\mathcal{K} \xrightarrow{\kappa} \mathbb{T}_{{\bf X}}\xrightarrow{\alpha} L[n] \quad \text{ in } \mathrm{Perf}(\bf X); \text{ and }$$ 
		\item the non-degeneracy condition: Locally on $\bf X$, where $L$ is trivialized, the induced $n$-shifted 1-form\footnote{We can locally identify the map $\alpha$ with the induced  1-form using the  trivialization of $L^{\vee}[n]$.} $\alpha: \mathbb{T}_{{\bf X}} \rightarrow \mathcal{O}_{{\bf X}}[n]$ is such that the map \begin{equation}
			d_{dR}\alpha|_{\mathcal{K}} \ \cdot:= \kappa^{\vee}[n] \circ( d_{dR}\alpha \ \cdot) \circ \kappa : \mathcal{K}\rightarrow\mathcal{K}^{\vee}[n]
		\end{equation} is a weak equivalence, where $d_{dR}\alpha \ \cdot$ is the induced contraction map $\mathbb{T}_{{\bf X}} \rightarrow \mathbb{L}_{{\bf X}}[n],$ and $\kappa^{\vee}: (\mathbb{T}_{{\bf X}})^{\vee}= \mathbb{L}_{{\bf X}} \rightarrow \mathcal{K}^{\vee}$ is the dual map. In that case, we say the $n$-shifted 2-form $d_{dR}\alpha$ is \bfem{non-degenerate on $\mathcal{K}$}.
	\end{itemize} Denote such a structure on $\bf X$ simply by $(\mathcal{K}, \kappa, L[n])$ or by $(L[n],\alpha)$, and call such local form a \emph{\textbf{$n$-contact form}}. Moreover, any such data satisfying all conditions above except the non-degeneracy is then called \bfem{$ n $-shifted pre-contact structure on $\bf X$.}

\end{definition}

From Proposition \ref{proposition_L as a complex of H^0 modules}, on a refined affine neighborhood, say $\spec A$ with $A$ a minimal standard form cdga, the perfect complexes $\mathbb{T}_A, \mathbb{L}_A$, when restricted to ${\spec H^0(A)}$, are both free finite complexes of $H^0(A)$-modules. In that case, Definition \ref{defn_shiftedcontact} will reduce to the following local descriptions, where  $\mathcal{K}$ is now just equivalent to the usual $\ker \alpha$ in $Mod_A$; and  $L$ in the splitting corresponds to the  line bundle generated by the Reeb vector field of the classical case.  More precisely, from \cite[$\S 3.2$]{kib}, when restricted to the (nice) local models,  we equivalently have the following proposition/definition. 
	\begin{proposition} \label{defn_contact strc on good affines}
		\emph{(Shifted contact structures with good affine models)} For a (minimal) standard form cdga $A$ and $n<0$, any $n$-shifted contact structure on  $\textbf{X}=\spec A$ can be strictified in the sense that the resulting contact data consists of \begin{itemize}
			\item a submodule $\mathcal{K}$ with the natural inclusion $i:\mathcal{K}\hookrightarrow Der(A)$
			such that $Cone(i)\simeq coker(i)$ is the quotient complex and of the form $L[n]$, with $L$ a line bundle; and
			\item  an  $n$-shifted 1-form $\alpha$  on $\spec A$  with the property that $\mathcal{K}\simeq \ker \alpha$ so that the $n$-shifted 2-form $d_{dR}\alpha$ is non-degenerate  on $\ker \alpha$.
		\end{itemize}Here $Der(A)=(\Omega^1_{A})^{\vee}=\Hom_{A}(\Omega^1_{A},A)$, where $ \Omega^1_{A} $ is the $A$-module K\"{a}hler differentials such that $ \Omega^1_{A}|_{\cspec H^0(A)} $ is represented as a (bounded) complex of free $H^0(A)$-modules. In that case,  one also has the following splitting over $ p\in \spec H^0(A): $	
	$$ 	Der(A)|_{\spec H^0(A)}=  \ker \alpha |_{\spec H^0(A)} \oplus L[n]|_{\spec H^0(A)}. $$ 
		
	\end{proposition} 
		Adopting the classical terminology, we sometimes call the sub-module $\ker \alpha$  above an \bfem{$ n $-shifted (strict) contact structure with the defining $ n $-contact form}  $ \alpha $.

In addition to the strictification result above, \cite{Berktav2024} provides several examples of contact derived stacks, outlined in the following theorem.
	
\begin{theorem}\label{ON EXAMPLES} \emph{(Examples of contact derived stacks)} Let $ \bf X $ be a derived Artin $\K$-stack locally of finite presentation. Denote by $ \mathbb{G}_a, \mathbb{G}_m $ the affine additive and multiplicative group schemes, respectively.
	\begin{enumerate}
		\item Let $\mathrm{T^*[n]} \bf X$ be the $n$-shifted cotangent stack. Then   the space $ \mathrm{J^1[n]{\bf X}:=T^*[n]}{\bf X }\times \mathbb{G}_a$\footnote{By abuse of notation, $ \mathbb{G}_a $ stands for the \emph{$n$-shifted affine line} corepresented by  the polynomial algebra on a variable in cohomological degree $ -n $. }, called the \emph{$n$-shifted 1-jet stack} of $\bf X$, carries an $n$-shifted contact structure.
		\item Let  $\pi_{{\bf X}}: \mathrm{T^*}{\bf X}\rightarrow \bf X$ be the natural projection. Given a prequantum 0-shifted Lagrangian fibration structure on $ \pi_{{\bf X}} $, there is a  $\mathbb{G}_m$-bundle on $ \mathrm{T^*}{\bf X}$ carrying a 0-shifted contact structure.
		\item Let  $\pi_{c_1(\mathcal{G})}: \mathrm{T^*}_{c_1(\mathcal{G})}\bf   X \rightarrow X$ be the $ c_1(\mathcal{G})$-twisted cotangent stack of $\bf X$, where $ c_1(\mathcal{G}) \in \mathcal{A}^1({\bf X}, 1) $ denotes the \emph{characteristic class} of a 0-gerbe $\mathcal{G}$ - a line bundle - on $\bf X$. Given a prequantum 0-shifted Lagrangian fibration structure on $ \pi_{c_1(\mathcal{G})} $, there is a  $\mathbb{G}_m$-bundle on $\mathrm{T^*}_{c_1(\mathcal{G})}\bf  X$ that carries a 0-shifted contact structure.
		\item Assume that $ G $ is a simple algebraic group
		over  $ \mathbb{K} $, and $ C $ be a  smooth and proper curve$/\mathbb{K} $. Consider the derived moduli stacks $LocSys_G(C), Bun_G(C)$ of \emph{flat $G$-connections on $C$}, \emph{principal $G$-bundles on $C$}. Then there is a  $\mathbb{G}_m$-bundle on $ LocSys_G(C) $ with a 0-shifted contact structure.
	\end{enumerate}
\end{theorem}
\paragraph{Local theory for shifted contact structures.} \label{section_darboux theorem}
It has been shown in \cite{kib}  that every $k$-shifted contact derived $\K$-scheme $\bf X$ is locally equivalent to $(\spec A, \alpha_0)$ for $A$ a minimal standard form cdga and $\alpha_0$ in Darboux form. In addition,  \cite{Berktav2024} extends this result to the case of derived Artin stacks. More precisely, we have the following result for derived schemes.

\begin{theorem} \cite[Thm. 3.12.]{kib}\label{contact darboux}
	Let $\bf X$ be  a $k$-shifted contact derived $\mathbb{K}$-scheme  for $k<0$, and $x\in \bf X$. Then there is a local contact model $\big(A,  \alpha_0 \big)$  and $p \in \spec H^0(A)$ such that $i: \spec A \hookrightarrow \bf X$ is an open inclusion with $i(p)=x$, \ $A$ is a standard form that is minimal at  $p$, and $\alpha_0$ is a $k$-shifted contact form on $\spec A$ such that $A, \alpha_0$ are in contact Darboux form.%, with $i^*(\alpha) \sim \alpha_0$ in the space of $k$-shifted 1-forms. 

\end{theorem} 

\begin{example}
When $k$ is odd, see \cite[Examaple 2.21-2.24]{Berktav2024}, one can find a minimal standard form cdga $A$, with $ `` $coordinates" $x^{-i}_j, y^{k+i}_j,z^k \in A$, and a Zariski open inclusion $f: \spec A \hookrightarrow \bf X$ so that
\begin{equation*}
	\alpha_0 = d_{dR}z^k +  \sum_{i,j} y_j^{k+i}d_{dR}x_j^{-i}
\end{equation*} defines a $k$-contact form on $\spec A$. Here, using a special element $H\in A^{k+1}$, called the \emph{Hamiltonian},  the \bfem{internal differential} $d$ on $A$ is given by  the equations \begin{align} 
d|_{A(0)}&=0; \quad dx_j^{-i} =  \dfrac{\partial H}{\partial y_j^{k+i}}; \quad dy_j^{k+i} =  \dfrac{\partial H}{\partial x_j^{-i}} \text{ for all } i,j; \text{ and } \nonumber \\ -kdz^k&= H+d\Big[\sum_{i,j} (-1)^{i} ix_j^{-i} y_j^{k+i} \Big].
\end{align}%\paragraph{Symplectization of a shifted contact derived scheme.} \label{section_symplectization}
{If $k=-1$}, in particular, we have  $A=A(0)[z^{-1}, y_1^{-1}, \dots, y_{m_0}^{-1}]$, with $A(0)$ a $\mathbb{K}$-algebra generated by $x_1^0,\dots, x_{m_0}^0$, such that $\vdim A= m_0 - (m_0+1)=-1$. Choosing an arbitrary Hamiltonian $H\in A(0)$, we let $dz^{-1}=H$ and $dy^{-1}_j= {\partial H}/{\partial x_j^{0}}$ and $dx^0_j=0 \ \forall j$.  Then,   a $(-1)$-contact form  can be written as $$\alpha_0=\dR z^{-1}+\sum_{1\leq j \leq m_0}y_j^{-1}\dR x_j^{0}.$$ In that case, $\ker \alpha_0$ is generated by the vector fields $ \partial/\partial y_j^{-1}$ { and } $y_j^{-1}\partial/\partial z^{-1} -   \partial/\partial x^{0}_j  \text{ for } 1\leq j\leq m_0.$ For more details, we refer to \cite{kib,Berktav2024}. 
\end{example}

\subsection{Overview of Legendrians in contact manifolds}

In this section, we provide some background from  contact geometry, leading to the concept of a Legendrian submanifold. For more details, we refer to standard sources like \cite{McDuff2017,Geiges}.

\begin{definition}
	Let $(M, \xi)$ be a contact manifold of dimension $2n+1$. A submanifold $L\subset M$ is called an \emph{\textbf{isotropic}} submanifold if $T_pL \subset \xi_p$ for all $p\in L.$ 
	
	Equivalently, we can formulate the isotropy condition for $L$ using locally defining contact forms as follows. \begin{definition} \emph{(Isotropic submanifolds - 2nd definition)}
		Let $\alpha$ be a local contact form for $\xi$, we say that a submanifold $i:L\hookrightarrow (M, \xi)$ is \emph{\textbf{isotropic}} if $i^*\alpha\equiv0.$ 
	\end{definition}
\end{definition}

\begin{observation}
Let $i:L\hookrightarrow (M, \xi)$ be an isotropic submanifold and  $\alpha$ be a  contact form (at least locally) defining $\xi$, then the condition for $L$ to be isotropic implies that $i^*(d_{dR}\alpha)=0$ as well. It means that $T_pL\subset\xi_p$ is an isotropic subspace of the $2n$-dimensional symplectic vector space $(\xi_p, d_{dR}\alpha|_{\xi_p})$. From linear algebra, we have $\dim T_pL \leq (\dim \xi_p)/2=n.$ This then leads to the following natural subobjects for which $\dim T_pL $ attains its maximum.
\end{observation}

\begin{definition}
	An isotropic submanifold $L\subset (M^{2n+1}, \xi)$ of maximal possible dimension $n$ is  called a \emph{\textbf{Legendrian}} submanifold. 
\end{definition}

In brief, Legendrians $L\subset (M, \xi)$ are the largest possible submanifolds whose tangent bundle is completely contained in the subbundle $ \xi $ of $TM$. In other words, they are the \textit{integral submanifolds} of $\xi$ of maximal dimension $ (\dim \xi_p)/2=n. $ In that case, $TL$ is a \emph{Lagrangian} subbundle of $i^*\xi.$
\begin{example} \label{Legendrians in 1-jet bundles}\emph{(Legendrians in 1-jet bundles)}
	Given a manifold $L$, define the \bfem{space $J^1(L)$ of 1-jets of germs of smooth functions} $f:L\rightarrow \R$ as follows: Denote by $\epsilon_p$ the space of germs at $p\in L,$ with the equivalence relation $f \sim g \iff f(p)=g(p) \text{ and } d_{dR}f_p=d_{dR}g_p.$  We write $j^1_pf$ for the equivalence class of a germ $f\in \epsilon_p.$ Then we set $J^1(L)= \bigcup_p \epsilon_p/\sim.$
	
	Now, using the mapping $j^1_pf \mapsto (f(p),d_{dR}f_p)$, we get an identification $J^1(L)\simeq \R \times T^*L.$ It follows that there exists a natural contact structure $\xi_{jet}$ on $J^1(L)$ given by $\xi_{jet}:= \ker (d_{dR}z-\lambda)$, where $\lambda$ is the Liouville form on $T^*L$, which is locally of the form $\lambda= \sum_i p^i d_{dR}q^i.$ 
\end{example}	
	With this example in hand, we have the following important observation, which leads to a universal neighborhood model for any Legendrian submanifold.
	\begin{observation}
		Given any manifold $L$, every smooth function $f\in C^{\infty}(L)$ gives a\emph{ Legendrian embedding} of $L$ into $(J^1(L), \xi_{jet})$ via the map $p \mapsto (f(p),d_{dR}f_p)$ (and the identification $J^1(L)\simeq \R \times T^*L$). Note in particular that the zero function $f\equiv0$ corresponds to the zero section $j^10 \simeq L \subset T^*L \subset T^*L \times \R \simeq J^1(L)$ of the 1-jet bundle $J^1(L).$ 
		
		Now, the following theorem ensures that for Legendrian submanifolds $L\subset (M, \xi)$, a neighborhood of the zero section of its 1-jet bundle $ J^1(L) $ provides a model for (a neighborhood of) $L$. Therefore, it follows that, as in the case of Lagrangian submanifolds, there are no local invariants for Legendrian submanifolds as well.
	\end{observation}

\begin{theorem}(Weinstein's Legendrian neighborhood theorem)
Every 	Legendrian submanifold $L\subset (M, \xi)$ has a neighborhood contactomorphic to a neighborhood of the zero section of its 1-jet bundle $ J^1(L) $.
\end{theorem}
In the upcoming sections, we will discuss derived versions of the above.

\subsection{Legendrians in contact derived  stacks}
In this section, we begin by formulating the concept of a  \emph{contact isotropic structure} on a morphism in the context of derived contact geometry. We then introduce  \emph{Legendrian structures}. Our definitions are essentially based on \cite[Definitions 1.1 \& 1.2.]{Calaque2016}.
 
\begin{definition} \label{defn_isotropic for Leg}
	Let $f: \bf Y \rightarrow {\bf X}$ be a morphism of derived Artin stacks with $\bf X$ carrying an $n$-shifted contact structure $(\mathcal{K}, \kappa, L[n])$.
	A \bfem{contact isotropic structure} on $f$ consists of \begin{itemize}
		\item a morphism $\rho: \mathbb{T}_{{\bf Y}} \rightarrow f^*(\mathcal{K})$ of perfect complexes commuting the diagram\footnote{Induced from the fiber sequence $\mathbb{T}_f \rightarrow \mathbb{T}_{{\bf Y}}\rightarrow f^*(\mathbb{T}_{\bf X})$, with $\mathbb{T}_f$ the \emph{relative/vertical tangent complex} of $f$.} %$ \mathbb{T}_{\bf Y} \rightarrow f^*(\mathbb{T}_{{\bf X}})$ %. That is, we have the homotopy commutative diagram	
\begin{equation} \label{defn_isotropic for leg diagram}
	\begin{tikzcd}
		\mathbb{T}_{{\bf Y}} \arrow[r, "\rho"] \arrow[rd] & f^*(\mathcal{K}) \arrow[d, "f^*\kappa"] \arrow[dl, phantom, "\lrcorner", very near start]\\
	\	& f^*(\mathbb{T}_{{\bf X}});                    
	\end{tikzcd}
\end{equation}
\item a homotopy between $0$ and $ f^*(\alpha) $ for any $n$-contact form $\alpha$.
	\end{itemize}

%Observe that, %locally on $\bf X$, where $L$ is trivial, 
%	for an  $n$-contact form $\alpha$, we have (locally) $ f^*(\alpha)\in \mathcal{A}^1({{\bf Y}},n)$. Since the image under $\rho$ lies in $ f^*(\mathcal{K})\simeq f^*(Cocone(\alpha))\simeq Cocone(f^*(\alpha))$, the diagram implies that $0\sim f^*(\alpha)$, where $f^*(\alpha): f^*(\mathbb{T}_{{\bf X}}) \rightarrow \mathcal{O}_{{\bf Y}}[n]$. I.e., %the pushforward of any ``horizontal" tangent vector in $ f^*(\mathbb{T}_{{\bf X}})$ under the map $f^*(\alpha): f^*(\mathbb{T}_{{\bf X}}) \rightarrow \mathcal{O}_{{\bf Y}}[n]$ will lie inside $\mathcal{K}\simeq Cocone(\alpha)$. 
%	 there is a path $\varLambda_{\alpha}$ between $0$ and $f^*(\alpha)$ such that $0\sim f^*(\alpha)$.
	
	We then  denote an \emph{isotropic structure} on $f: \bf Y \rightarrow {\bf X}$ equivalently by $(\rho, \varLambda)$, where $\varLambda$ is the map $$\mathcal{A}^1_{loc}({\bf X}, n)\longrightarrow \bigsqcup_{\alpha' \in \mathcal{A}^1_{loc}({\bf X}, n)}Path_{\mathcal{A}^{1}_{loc}({\bf Y}, n)} {\big(0, f^*(\alpha')\big), \ \alpha \mapsto \varLambda_{\alpha}}.$$% homotopy $\varLambda_{{\bf Y}}$ from $0$ to the composition $f^*(\mathcal{K})\rightarrow f^*(\mathbb{T}_{{\bf X}})\rightarrow \mathcal{O}_{{\bf Y}}[k]$. 
\end{definition}

\begin{observation}
Let $(\rho, \varLambda)$ be an isotropic structure   on $f: {\bf Y}\rightarrow {\bf X}$ with an $n$-shifted contact target $({\bf X}; \mathcal{K}, \kappa, L[n])$. Choose a locally defining $n$-contact form $\alpha$ and a path $\varLambda_{\alpha}$ from $0$  to $f^*(\alpha)$. By definition, we require the morphism $\rho: \mathbb{T}_{{\bf Y}} \rightarrow f^*(\mathcal{K})$ to make
Diagram \ref{defn_isotropic for leg diagram}
commute. Then locally on $\bf X$, the path $\varLambda_{\alpha}$ gives a homotopy between 0 and the composition  	\begin{equation} \label{induced homotopy1}
	\mathbb{T}_{{\bf Y}}\rightarrow f^*(\mathcal{K})\rightarrow f^*(\mathbb{T}_{{\bf X}})\xrightarrow{f^*(\alpha)} \mathcal{O}_{{\bf Y}}[n].
	\end{equation}
\end{observation}

Now, using the identity $d_{dR}\circ f^* \simeq f^* \circ d_{dR}$ and the path $\varLambda_{\alpha}$, we obtain a path $d_{dR}\varLambda_{\alpha}$ between 0 and $f^*(d_{dR}\alpha): f^* (\mathbb{T}_{\bf X}) \wedge f^*(\mathbb{T}_{\bf X}) \rightarrow \mathcal{O}_{{\bf Y}}[n]$. Then using the  triangle $\mathbb{T}_f \rightarrow \mathbb{T}_{{\bf Y}}\rightarrow f^*(\mathbb{T}_{\bf X})$, this path (locally) induces a  homotopy from 0 to the map $\mathbb{T}_{\bf Y} \wedge \mathbb{T}_{\bf Y}  \rightarrow \mathcal{O}_{{\bf Y}}[n] $ defined as the composition $$\mathbb{T}_{\bf Y} \wedge \mathbb{T}_{\bf Y} \rightarrow f^* (\mathbb{T}_{\bf X}) \wedge f^*(\mathbb{T}_{\bf X}) \rightarrow \mathcal{O}_{{\bf Y}}[n].$$  %Using the triangle above one more time, the path $d_{dR}\varLambda_{\alpha}$ also induces a  local homotopy between 0 and the composition $$\mathbb{T}_f \otimes \mathbb{T}_{\bf Y}\rightarrow \mathbb{T}_{\bf Y} \wedge \mathbb{T}_{\bf Y} \rightarrow f^* (\mathbb{T}_{\bf X}) \wedge f^*(\mathbb{T}_{\bf X}) \rightarrow \mathcal{O}_{{\bf Y}}[k].$$ %Recall that the map $\mathbb{T}_f \rightarrow f^*(\mathbb{T}_{\bf X})$ comes with a natural homotopy to $0$, and hence we get another homotopy between 0 and $ \mathbb{T}_f \otimes \mathbb{T}_{\bf Y}\rightarrow \mathcal{O}_{{\bf Y}}[k] $. Combining all of these homotopies, we obtain a loop at 0 in the space $Map_{L_{QCoh({\bf Y})}} (\mathbb{T}_f \otimes \mathbb{T}_{\bf Y}, \mathcal{O}_{{\bf Y}}[k])$. That is, we get an element in $\pi_1(Map_{L_{QCoh({\bf Y})}} (\mathbb{T}_f \otimes \mathbb{T}_{\bf Y}, \mathcal{O}_{{\bf Y}}[k]); 0)\simeq [\mathbb{T}_f \otimes \mathbb{T}_{\bf Y}, \mathcal{O}_{{\bf Y}}[k-1]]$. Thus, using adjunction, 
Now, from Observation \ref{observation_gettting LAG non-deg } (with $ d_{dR}\alpha $ in place of $\omega_{\bf X}^0$)\footnote{Keep in mind that $d_{dR}\alpha \cdot: \mathbb{T}_{{\bf X}} \rightarrow \mathbb{L}_{{\bf X}}[n]$ is a quasi-isomorphism on $\mathcal{K}.$}, we locally get the induced morphism of perfect complexes on ${\bf Y}$ (depending on $\varLambda_{\alpha}$)
\begin{equation} 
\chi_{\varLambda_{\alpha}}: \mathbb{T}_f \rightarrow \mathbb{L}_{{\bf Y}}[n-1].
\end{equation}

 This observation leads to the formulation of a \textit{non-degeneracy condition} for  (contact) isotropic structures given in  Definition \ref{defn_isotropic for Leg} and  the notion of a \emph{Legendrian structure}. In that respect, we introduce the following definition.

%\[\begin{tikzcd}
%\	& f^*(\mathcal{K}) \arrow[d]  \arrow[rd, dashed] &                   \\
%	\mathbb{T}_{{\bf Y}} \arrow[ru]   \arrow[r]   & f^*(\mathbb{T}_{{\bf X}}) \arrow[ul, phantom, "\lrcorner", very near start] \arrow[r]                   & {\mathbb{T}_f[1]},
	%\arrow[dr, phantom, "\ulcorner", very near start]
%\end{tikzcd}
% \]

\begin{definition} \label{defn_LEG strc}
Let $f: \bf Y \rightarrow \bf X$ be a morphism of derived Artin stacks with an $n$-shifted contact target  $({\bf X}; \mathcal{K}, \kappa, L[n])$. A \bfem{Legendrian structure} on $f$ consists of \begin{itemize}
	\item a contact isotropic structure $\rho: \mathbb{T}_{{\bf Y}} \rightarrow f^*(\mathcal{K})$ on $f$; and
	\item  a fiber sequence\footnote{Recall from \cite{Calaque2016}, for a generic map $\phi: U \rightarrow V$ of perfect complexes with a degree $n$ (skew-symmetric) pairing on $V$, we have the induced null-homotopic sequence $U\rightarrow V \rightarrow U^{\vee}[n]$. In our case with $\phi=\rho$, $V:=f^*(\mathcal{K})$ comes with a pairing $d_{dR}\alpha|_{\mathcal{K}} \cdot$ induced by an $n$-contact form $\alpha$.} $\mathbb{T}_{{\bf Y}} \xrightarrow{\rho} f^*(\mathcal{K}) \rightarrow \mathbb{L}_{{\bf Y}}[n]$ \ (the \emph{non-degeneracy condition} for $\rho$). 
	%$f^*(\mathcal{K}) \rightarrow f^*(\mathbb{T}_{{\bf X}}) \rightarrow f^*(\mathcal{K})^{\vee}[n]$ \ (the \emph{non-degeneracy condition} for $\rho$). 
\end{itemize} We then say $\bf Y$ is a \textit{Legendrian} in $({\bf X}; \mathcal{K}, \kappa, L)$.
%with the data of a commuting square
%\begin{equation} 
%\begin{tikzpicture}
%\matrix (m) [matrix of math nodes,row sep=1.5 em,column sep=3 em,minimum width=1 em] {
%	\mathbb{T}_f[1]  &  	\mathbb{L}_{{\bf Y}} [k] \\
%	  f^*(\mathcal{K}) & (f^*(\mathcal{K}))^{\vee}[k].\\};
%\path[-stealth]
%(m-1-1) edge  node [left] { } (m-1-2)
%edge  node [above] { } (m-2-1)

%(m-2-1) edge  node [below] { } (m-2-2)
%(m-2-2) edge  node [below] { } (m-1-2);

%\end{tikzpicture}
%\end{equation} Here the vertical maps are obtained from the (dual of) morphism $\rho$ and the triangle $\mathbb{T}_f \rightarrow \mathbb{T}_{{\bf Y}}\rightarrow f^*(\mathbb{T}_{\bf X})$; and the horizontal maps are $\chi_{\varLambda_{\alpha}}[1]$ and the induced map $f^*(d_{dR}\alpha)\cdot$.
\end{definition}

In other words, a \bfem{Legendrian structure} on such $f$ is defined to be a \emph{non-degenerate isotropic structure on the morphism $\rho: \mathbb{T}_{{\bf Y}} \rightarrow f^*(\mathcal{K})$} in the sense of \cite[Definitoin 1.2]{Calaque2016}.
\begin{observation}\label{observation_non-deg isot str induces weak eqv on pullback}
	Suppose that we have a morphism $f: {\bf Y} \rightarrow ({{\bf X}}; \mathcal{K}, \kappa, L[n])$ of derived Artin stacks, with a shifted contact target, and  a contact isotropic structure $\rho: \mathbb{T}_{{\bf Y}} \rightarrow f^*(\mathcal{K})$ on $f$. By definition, we have the induced morphism $f^*(d_{dR}\alpha|_{\mathcal{K}}): \wedge^2f^*(\mathcal{K}) \rightarrow \mathcal{O}_{{\bf Y}}[n]$ and the pairing $ f^*(d_{dR}\alpha|_{\mathcal{K}}) \ \cdot : f^*(\mathcal{K})\rightarrow f^*(\mathcal{K}^{\vee}[n])$, not necessarily non-degenerate. 
	
	If, in addition,  the contact isotropic structure $\rho$ is non-degenerate in the sense of Definition \ref{defn_LEG strc}, the induced null-homotopic sequence $\mathbb{T}_{{\bf Y}} \xrightarrow{\rho} f^*(\mathcal{K}) \rightarrow \mathbb{L}_{{\bf Y}}[n]$ and its shifted dual are two fiber sequences fitting into the diagram
	
	\begin{equation}
		\begin{tikzcd}
			\mathbb{T}_{{\bf Y}} \arrow[d, equal] \arrow[r] & f^*(\mathcal{K}) \arrow[r] \arrow[d, "f^*(d_{dR}\alpha|_{\mathcal{K}}) \ \cdot", dashed] & {\mathbb{L}_{{\bf Y}}[n]} \arrow[d, equal] \\
			\mathbb{T}_{{\bf Y}} \arrow[r]                                & { f^*(\mathcal{K}^{\vee}[n])} \arrow[r]                                                  & {\mathbb{L}_{{\bf Y}}[n].}                               
		\end{tikzcd}
	\end{equation}From \cite[Lemma 1.3]{Calaque2016}, the map $ f^*(d_{dR}\alpha|_{\mathcal{K}}) \ \cdot : f^*(\mathcal{K})\rightarrow f^*(\mathcal{K}^{\vee}[n])$ is then a weak equivalence (i.e. a quasi-isomorphism), hence a non-degenerate pairing.
\end{observation}
\begin{observation}
	\emph{{(An alternative way of stating the non-degeneracy.)}}
	%The non-degeneracy condition above can be equivalently stated in terms of $n$-contact forms as follows: 
	Let $f$ be as in Definition \ref{defn_LEG strc} and $\rho$ an isotropic structure on $f$. Then the map $\mathbb{T}_{{\bf Y}} \rightarrow f^*(\mathbb{T}_{{\bf X}})$	factorizes as $ \mathbb{T}_{{\bf Y}}\rightarrow f^*(\mathcal{K})\rightarrow f^*(\mathbb{T}_{{\bf X}})$ up to homotopy. Let $\mathbb{T}_{\rho}:=\mathrm{hofib}\left(\rho: \mathbb{T}_{{\bf Y}}\rightarrow f^*(\mathcal{K})\right)$ and consider the fiber sequence
	\begin{equation}\label{hofib of rho}
		\mathbb{T}_{\rho}\rightarrow \mathbb{T}_{{\bf Y}}\xrightarrow{\rho} f^*(\mathcal{K}).
	\end{equation} Notice that if $\mathbb{T}_{{\bf Y}} \xrightarrow{\rho} f^*(\mathcal{K}) \rightarrow \mathbb{L}_{{\bf Y}}[n]$ is also a fiber sequence, then one obtains $\mathbb{T}_{\rho}[1]\simeq \mathbb{L}_{{\bf Y}}[n]$ using the sequence in (\ref{hofib of rho}), and vice versa. So, we may impose the equivalence \begin{equation}
	\boxed{\mathbb{T}_{\rho}\xrightarrow{\sim} \mathbb{L}_{{\bf Y}}[n-1]}
\end{equation} as an \bfem{alternative condition for the non-degeneracy of $\rho$}, leading to the following definition. 
\end{observation} 
\begin{definition}
Let $f: \bf Y \rightarrow \bf X$ be a morphism of derived Artin stacks with an $n$-shifted contact target  $({\bf X}; \mathcal{K}, \kappa, L[n])$ and $\rho: \mathbb{T}_{{\bf Y}}\rightarrow f^*(\mathcal{K})$ an isotropic structure on $f$. We say $\rho$ is  \bfem{non-degenerate} if the induced morphism $ \chi_{\rho}:\mathbb{T}_{\rho} \rightarrow \mathbb{L}_{{\bf Y}}[n-1]$ is a weak equivalence.

In that case, we may also refer to that condition by simply saying the \bfem{map $ \chi_{\varLambda_{\alpha}}: \mathbb{T}_f \rightarrow \mathbb{L}_{{\bf Y}}[n-1] $ is a weak equivalence or non-degenerate on $\mathbb{T}_{\rho}$.}
\end{definition}

	Let us now clarify the existence of a canonical map $\mathbb{T}_{\rho}\rightarrow \mathbb{L}_{{\bf Y}}[n-1]$. Recall that once we choose an $n$-contact form $\alpha$ defining $\mathcal{K}$, we locally get the induced morphism 
	$ 	\chi_{\varLambda_{\alpha}}: \mathbb{T}_f \rightarrow \mathbb{L}_{{\bf Y}}[n-1] $ of perfect complexes on ${\bf Y}$. Combining the universality of $\mathbb{T}_f$ with Diagram \ref{defn_isotropic for leg diagram}, we obtain the canonical morphism $ \mathbb{T}_{\rho} \rightarrow \mathbb{T}_{f} $ commuting diagram
	\begin{equation}
		\begin{tikzcd}
			\mathbb{T}_{\rho} \arrow[r, dashed] \arrow[rd] & \mathbb{T}_f \arrow[rd, "\lrcorner", phantom, near start] \arrow[d, "\pi_1"'] \arrow[r] & 0 \arrow[d]                              \\
			& \mathbb{T}_{{\bf Y}} \arrow[rd, "\rho"'] \arrow[r, "\pi_2"]                             & f^*(\mathbb{T}_{{\bf X}})                \\
			&                                                                                         & f^*(\mathcal{K}) \arrow[u, "f^*\kappa"'],
		\end{tikzcd}
	\end{equation}and hence the composition  $ \mathbb{T}_{\rho} \rightarrow \mathbb{T}_{f} \xrightarrow{\chi_{\varLambda_{\alpha}}} \mathbb{L}_{{\bf Y}}[n-1]$.

 \begin{observation}\label{observ_point with can contact data} %Let $ {\bf X}=\spec \mathbb{K}. $
 	%A Legendrian structure on $f: {\bf Y} \rightarrow \star_{n}$ gives a $(n-1)$-shifted contact structure on $\bf Y.$ 	Here, 
 	Denote by $\star_{n}$  the point $ {\bf X}:=\spec \mathbb{K} $ equipped with its canonical $n$-shifted contact structure, where the $n$-contact data $(\mathcal{K}, \kappa, L[n])$ for $\star_{n}$ is defined by the  $n$-shifted 1-form $\alpha \equiv 0.$ More precisely, the data consists of 
 	\begin{itemize}
 		\item the perfect complex $\mathcal{K}:= Cocone(\alpha)\simeq\mathbb{T}_{\bf X}=0$; and \item the trivial fiber sequence
 		$ Cocone(\alpha) \xrightarrow{\kappa} \mathbb{T}_{{\bf X}} \rightarrow 0{[n]}=:L[n] $ of perfect complexes,
 	where $ \kappa $ is  the zero map and %$Cone(\kappa)=\mathcal{O}_{{\bf X}}[n]\simeq \mathbb{K}[n],$ 
 	such that %the map $\kappa=0: \mathbb{T}_{\bf X} \rightarrow \mathbb{T}_{\bf X}$; and $L=0$, along with  % Notice that $Cone(\kappa)=\mathbb{T}_{\bf X}[-1] \oplus \mathbb{T}_{\bf X}$ is trivial and quasi-isomorphic to $\mathcal{O}_{\bf X}[k]$ for $k<0$ .
 the non-degeneracy of $d_{dR}\alpha$ on $\mathcal{K}$ trivially holds.
 	\end{itemize}We denote by such canonical structure $ (\star_{n}; \alpha \equiv 0, 0[n]). $
 
 Now, we can prove the following result.
 \end{observation}

\begin{proposition}\label{prop_Leg on f with point target}
Let $(\star_{n}; \alpha \equiv 0, 0[n])$ be as above. Any Legendrian structure $(\rho,\varLambda_{\alpha})$ on $f: {\bf Y} \rightarrow \star_{n}$ induces a natural $(n-1)$-shifted exact symplectic structure on $\bf Y$ and an $(n-1)$-shifted contact structure on ${\bf Y} \times \mathbb{A}^1[n-1]$, where $\mathbb{A}^1[n-1]$ is the \emph{$(n-1)$-shifted affine line} corepresented by the polynomial algebra on a variable in cohomological degree $ 1-n $.

%If, in addition, $Cocone(\varLambda_{\alpha})\simeq 0$, then the induced pre-contact data of $\bf Y$ can be promoted to an $(n-1)$-shifted contact structure on $\bf Y$.  %where $\star_{k}$ comes with its canonical $k$-shifted contact structure as above.
\end{proposition}
 \pf 
 
 Since $\star_{n}$ is equipped with the canonical $n$-shifted contact structure as above, any isotropic structure on $f$ is of the form $(\rho=0: \mathbb{T}_{{\bf Y}} \rightarrow f^*(\mathcal{K})=0, \varLambda)$, where $\varLambda_{\alpha}: 0 \rightsquigarrow f^*\alpha.$
 
 In our  case, $\alpha \equiv 0$ implies that $\varLambda_{\alpha}$ is a loop at 0 in the space $\mathcal{A}^1({\bf Y}, n)$. Therefore, it  defines an element in $\pi_1(\mathcal{A}^1({\bf Y}, n))\simeq \pi_0(\mathcal{A}^1({\bf Y}, n-1))$. Denote the corresponding 1-form of degree $(n-1)$ also by $\varLambda_{\alpha}: \mathbb{T}_{{\bf Y}}\rightarrow \mathcal{O}_{{\bf Y}}[n-1]$. Likewise, $d_{dR}\varLambda_{\alpha}$ is a loop at 0 due to $f^*(d_{dR}\alpha)\sim 0.$ Thus, it  gives an element in $ \pi_0(\mathcal{A}^2({\bf Y}, n-1))$, hence an $(n-1)$-shifted exact pre-symplectic structure on $\bf Y$.
 
 As the operator $d_{dR}$  commutes with the identifications $\pi_0(\mathcal{A}^p({\bf Y}, n-1))\simeq \pi_1(\mathcal{A}^p({\bf Y}, n))$,  the induced map $ \chi_{\varLambda_{\alpha}}$ is nothing else than the contraction map $d_{dR}\varLambda_{\alpha} \cdot: \mathbb{T}_f\simeq\mathbb{T}_{{\bf Y}} \rightarrow \mathbb{L}_{{\bf Y}}[n-1]$, which is a weak equivalence on $ \mathbb{T}_{\rho} $ due to the Legendrian structure  on $f$. 
 
 Since $\rho: \mathbb{T}_{{\bf Y}} \rightarrow f^*(\mathcal{K})=0$, we get $ \mathbb{T}_{{\rho}}=\mathrm{hofib}(\rho)=\mathbb{T}_{{\bf Y}}.$ Thus, the contraction morphsim $d_{dR}\varLambda_{\alpha} \cdot: \mathbb{T}_{{\bf Y}} \rightarrow \mathbb{L}_{{\bf Y}}[n-1]$ happens to be non-degenerate on $ \mathbb{T}_{{\bf Y}} $, implying that $d_{dR}\varLambda_{\alpha}$ defines an  $(n-1)$-shifted (exact) symplectic structure on $\bf Y$.

Now, let us consider the product space  ${\bf Y} \times \mathbb{A}^1[n-1]$ defined as the pullback diagram
\[\begin{tikzpicture}
	\matrix (m) [matrix of math nodes,row sep=3em,column sep=5 em,minimum width=2 em] {
		\mathcal{Z}:={\bf Y} \times \mathbb{A}^1[n-1]	   & \mathbb{A}^1[n-1]  \\
	 ({\bf Y}, d_{dR}\varLambda_{\alpha}) &  (\star_{n}; \alpha \equiv 0, 0[n]), \\};
	\path[-stealth]
	(m-1-1) edge  node [left] { $ pr_1 $} (m-2-1)
	(m-1-1) edge  node [above] { $ pr_2 $} (m-1-2)
	%(m-1-1) edge  node [below] {} node [below] {{\small stacks}} (m-2-2)
	%(m-1-1) edge  node [below] {} node [below] {{\small  higher stacks}} (m-3-2)
	(m-2-1) edge  node  [above] {$  $ } (m-2-2)
	
	(m-1-2) edge  node [right] {$  $ } (m-2-2);
	%edge [dashed,-] (m-2-1);
\end{tikzpicture}\] where $\bf Y$ carries the   $(n-1)$-shifted exact symplectic structure $d_{dR}\varLambda_{\alpha}$, and $\star_{n}$ is equipped with the canonical $n$-shifted contact structure as above. The pullback diagram implies that we have the fiber-cofiber sequence \[
pr_1^*\mathbb{T}_{{\bf Y}} \rightarrow \mathbb{T}_{\mathcal{Z} }\rightarrow pr_2^* \mathbb{T}_{\mathbb{A}^1[n]},\] where $ \mathbb{T}_{\mathbb{A}^1[n-1]} $ is an invertible quasi-coherent sheaf, and hence,  $pr_2^* \mathbb{T}_{\mathbb{A}^1[n]} \simeq L [n]$ with $L$  a line bundle, leading to an $(n-1)$-shifted pre-contact structure.

Letting $\sigma:= \dR t + \Lambda_{\alpha}$, where $t$ is the degree $-n+1$ variable for the algebra corepresenting $ \mathbb{A}^1[n-1]$, one can conclude from the construction that $\dR \sigma$ is non-degenerate on $ pr_1^*\mathbb{T}_{{\bf Y}}$ and $Cocone(\sigma) \simeq pr_1^*\mathbb{T}_{{\bf Y}}$ (cf. \cite[Lemma 4.1]{Berktav2024}); and hence, $\sigma$ defines an $(n-1)$-shifted contact structure on the product $\mathcal{Z}={\bf Y} \times \mathbb{A}^1[n-1]$. 
\epf

\subsection{Example: The zero section morphism of 1-jet stacks} In this section, we show that the zero section morphism of a 1-jet stack carries a natural Legendrian structure, providing the derived extension of  Example \ref{Legendrians in 1-jet bundles}. 

Let $  T^*[n]X $ be the $n$-shifted cotangent stack of the derived Artin stack $ X$ (locally of finite presentation). Recall from \cite{Berktav2024} that the \emph{\textbf{$n$-shifted 1-jet stack of $X$}} is described  by the pullback diagram \begin{equation} \label{diagram_cartesian spaces as pullbacs}
	\begin{tikzpicture}
		\matrix (m) [matrix of math nodes,row sep=3em,column sep=5 em,minimum width=2 em] {
			J^1[n]X:=T^*[n]X \times \mathbb{A}^1[n]	   & \mathbb{A}^1[n]  \\
			T^*[n]X   &  * \\};
		\path[-stealth]
		(m-1-1) edge  node [left] { $ pr_1 $} (m-2-1)
		(m-1-1) edge  node [above] { $ pr_2 $} (m-1-2)
		%(m-1-1) edge  node [below] {} node [below] {{\small stacks}} (m-2-2)
		%(m-1-1) edge  node [below] {} node [below] {{\small  higher stacks}} (m-3-2)
		(m-2-1) edge  node  [above] {$  $ } (m-2-2)
		
		(m-1-2) edge  node [right] {$  $ } (m-2-2);
		%edge [dashed,-] (m-2-1);
	\end{tikzpicture}
\end{equation}Here, $\mathbb{A}^1[n]$ is the \emph{$n$-shifted affine line} corepresented by the polynomial algebra on a variable in cohomological degree $ -n $,  with the  projection map being the composition $J^1[n]X\rightarrow T^*[n]X\rightarrow X.$  Moreover, from Diagram \ref{diagram_cartesian spaces as pullbacs}, we have the identifications
\begin{equation}\label{identifications by pullback}
\mathbb{L}_{J^1[n]X} \simeq pr_1^*\mathbb{L}_{T^*[n]X} \oplus pr_2^* \mathbb{L}_{\mathbb{A}^1[n]} \text{ and } \mathbb{T}_{J^1[n]X} \simeq pr_1^*\mathbb{T}_{T^*[n]X} \oplus pr_2^* \mathbb{T}_{\mathbb{A}^1[n]} .
\end{equation}
Equivalently, we have the exact triangle \begin{equation}\label{triangle giving precontact}
pr_1^*\mathbb{T}_{T^*[n]X} \rightarrow \mathbb{T}_{J^1[n]X} \rightarrow pr_2^* \mathbb{T}_{\mathbb{A}^1[n]}.
\end{equation} Notice that from the triangle (\ref{triangle giving precontact}), we have a natural morphism $pr_1^*\mathbb{T}_{T^*[n]X} \rightarrow \mathbb{T}_{J^1[n]X}$ with the cofiber $ pr_2^* \mathbb{T}_{\mathbb{A}^1[n]}. $ Since $ \mathbb{T}_{\mathbb{A}^1[n]} $ is an invertible quasi-coherent sheaf, we get $pr_2^* \mathbb{T}_{\mathbb{A}^1[n]} \simeq L [n]$, where $L$ is a line bundle.

Recall that using the observations above, we can equip  $ J^1[n]{X}$ with a canonical $n$-shifted contact structure. More precisely, we have:\begin{lemma} \label{lemma_can contact for 1-jet} \cite[Thm. 4.2]{Berktav2024}
	 The \bfem{canonical $n$-shifted contact data for $ J^1[n]{X}$} consists of \begin{itemize}
	\item the perfect complex $\mathcal{K}:= pr_1^*\mathbb{T}_{T^*[n]X}$ and the  morphism $\kappa: pr_1^*\mathbb{T}_{T^*[n]X} \rightarrow \mathbb{T}_{J^1[n]X}$ in (\ref{triangle giving precontact});
	\item the cofiber $ pr_2^* \mathbb{T}_{\mathbb{A}^1[n]}$  in (\ref{triangle giving precontact}) as the cone of $\kappa$, where  we  have an equivalence  $pr_2^* \mathbb{T}_{\mathbb{A}^1[n]} \simeq L [n]$ with $L$ a line bundle; and
\item  the equivalence $\mathcal{K}\simeq Cocone(\alpha)$ by \cite[Lemma 4.1]{Berktav2024}, with   the global $n$-contact form $\alpha$ given as
\begin{equation}
	\alpha:= -d_{dR}z + \lambda \in \Gamma\big(J^1[n]X, \mathbb{L}_{J^1[n]X}[n] \big),
\end{equation} where $\lambda$ is the \emph{Liouville 1-form} on $T^*[n]X$, and we simply write $z, \lambda$ instead of $pr_2^*z, pr_1^*\lambda,$ respectively. For the precise definition of $\lambda$, see \cite{Calaque2016}.
\end{itemize}
\end{lemma}
Now, we prove the following result.
\begin{theorem}\label{thm_on zero section}
	Let $ J^1[n]{X}$ be the \emph{$n$-shifted 1-jet stack} of $X$ with the canonical $n$-shifted contact structure as in Lemma \ref{lemma_can contact for 1-jet}. Then the zero section morphism $j: X \rightarrow {J^1[n]}{ X}$ has a natural  Legendrian structure.
\end{theorem}
\pf
Using the natural \emph{$n$-shifted contact data} on $ {J^1[n]}{ X} $  in Lemma \ref{lemma_can contact for 1-jet},   we observe:
\begin{itemize}
	\item \emph{There is an isotropic structure $\rho$ on $j$:} It should be noted that the map $j$ is the composition \begin{equation}
		j=\iota \circ \iota_{{X}}: X\xrightarrow{\iota_{{X}}} T^*[n]X \xrightarrow{\iota} J^1[n]X,
	\end{equation} where both $ \iota_{{X}}, \iota $ are the corresponding zero section morphisms. 

Consider the pullback  $j^*\kappa: j^*(\mathcal{K}) \rightarrow j^*(\mathbb{T}_{J^1[n]X})$. Since $pr_1\circ \iota \sim id$ and $\mathcal{K}= pr_1^*\mathbb{T}_{T^*[n]X}$, we get $$ j^*(\mathcal{K})=j^*(pr_1^*\mathbb{T}_{T^*[n]X})\simeq  (pr_1\circ j)^*\mathbb{T}_{T^*[n]X} \simeq (pr_1\circ \iota \circ \iota_{{X}})^*\mathbb{T}_{T^*[n]X} \simeq  \iota_{{X}}^*\mathbb{T}_{T^*[n]X}. $$  Thus, we have a homotopy between $ j^*\kappa $ and the map $\iota_{{X}}^*\mathbb{T}_{T^*[n]X}\rightarrow j^*\mathbb{T}_{J^1[n]X},$ written $ j^*\kappa \sim (\iota_{{X}}^*\mathbb{T}_{T^*[n]X}\rightarrow j^*\mathbb{T}_{J^1[n]X})$.

Choose $\rho$ to be natural map $\mathbb{T}_X \rightarrow  \iota_{{X}}^*\mathbb{T}_{T^*[n]X}$, commuting the diagram 
\begin{equation} 
	\begin{tikzcd}
		\mathbb{T}_{{X}} \arrow[r] \arrow[rd] & {j^*(\mathcal{K})\simeq \iota_{{X}}^*\mathbb{T}_{T^*[n]X}} \arrow[dl, phantom, "\lrcorner", very near start] \arrow[d] \\
	\	& {j^*(\mathbb{T}_{J^1[n]X})}.                                                           
	\end{tikzcd}
\end{equation}
Observe also that using the global $n$-contact form $\alpha$ given in Lemma \ref{lemma_can contact for 1-jet}, we get \begin{equation*}
	j^*\alpha \sim (\iota \circ \iota_{{X}})^* \alpha \sim \iota_{{X}}^* (\iota^* \alpha) \sim 0.
\end{equation*} This is because $\iota^* \alpha$ is an $1$-form of degree $n$ on $T^*[n]X$ and the zero section morphism $\iota_{{X}}$ corresponds to the data of the identity map $id_X$ together with the zero section of $\mathbb{L}_X [n] $. Tautologically, this means that the pull-back  $\iota^* \alpha$ is nullhomotopic.

\item \emph{The isotropic structure is non-degenerate:} %We need to show that $\mathbb{T}_j  \rightarrow \mathbb{L}_X[n-1]$ is an equivalence on $j^*(\mathcal{K})[-1]$. Equivalently, i
It suffices to show that the induced null-homotopy sequence  
	$ \mathbb{T}_X \rightarrow j^*(\mathcal{K}) \rightarrow \mathbb{L}_X[n] $
is a fiber sequence.

We already showed that there is an identification $ j^*(\mathcal{K})=j^*(pr_1^*\mathbb{T}_{T^*[n]X}) \simeq \iota_{{X}}^*\mathbb{T}_{T^*[n]X},$ with the natural map  $\iota_{{X}}^*\mathbb{T}_{T^*[n]X}\rightarrow j^*\mathbb{T}_{J^1[n]X}.$ Note also that from \cite[Theorem 2.2]{Calaque2016}, the isotropic structure on $\iota_{{X}}$ is non-degenerate, and hence we  have a fiber sequence
$$ \mathbb{T}_X\rightarrow \iota_{{X}}^*\mathbb{T}_{T^*[n]X} \rightarrow \mathbb{L}_X[n], $$ 
which gives the desired result due to the identification $ j^*(\mathcal{K}) \simeq \iota_{{X}}^*\mathbb{T}_{T^*[n]X}.$
\end{itemize}\epf

\section{Local models for Legendrians in contact derived schemes}\label{section_local models and LEG darboux}
This section first presents {affine models}, called \bfem{Legendrian-Darboux forms} (cf. Construction \ref{Leg model}), for the Legendrians in  shifted contact derived $ \mathbb{K}$-schemes, and then provides a \bfem{Legendrian neighborhood theorem} (cf. Theorem \ref{Legendrian-Darboux}) %We first introduce explicit local models, called \bfem{Legendrian-Darboux forms} (cf. Construction \ref{Leg model}), %for the case when $k<0$ is odd with $k\equiv 1\mod 4$. %cases $(1) \ k<0 \text{ with } k\not\equiv 3\mod 4$ and $(2) \ k<0 \text{ with } k\equiv 3\mod 4$. 
%then give the proof of Theorem \ref{THM3}
, which essentially says that \emph{any such structure can be locally modeled on  the Legendrian Darboux forms}.

We should  note that we will provide the proofs \emph{only for} negative odd shifts satisfying $k\equiv 1\mod 4$% the case $(1)$ above
, and give a brief sketch for the other cases, $ k\equiv 3\mod 4 $ or $k$ even. %is a variation on \cite[Theorem 3.7]{JS}. Additionally,  our prototype model, Construction \ref{Leg model}, will be  based on \cite[Example 3.3]{JS}. 
Additionally,  our constructions will be based on \cite[Example 3.3 \& Theorem 3.7]{JS}.    

\subsection{Legendrian-Darboux forms} %Recall that given  a locally finitely presented  derived Artin stack $\bf X$ with an $n$-shifted contact data $(\mathcal{K}, \kappa, L[n])$, the main definitions with nice local models can be strictified as in Proposition \ref{defn_contact strc on good affines}. In what follows, we assume such strict affine models. %reduce to the case, where  $\mathcal{K}$ is locally just the usual $\ker \alpha$ in $Mod_A$ with the map $\ker \alpha \hookrightarrow \mathbb{T}_A$ of complexes of $H^0(A)$-modules for $A$ a standard from cdga; and  $L$  corresponds to the  line bundle generated by the Reeb vector field of the classical case.  That is, we have:	

%\begin{definition}\emph{(Affine case with standard forms, see Proposition \ref{defn_contact strc on good affines})}
% An \bfem{$ n $-shifted  contact structure on $\textbf{X}=\spec A$}, with $A$ a standard form cdga, is a submodule $\mathcal{K}$ of $\mathbb{T}_A$ such that $\mathcal{K}\simeq \ker \alpha$ for a $n$-shifted 1-form $\alpha$    with the property that the 2-form $d_{dR}\alpha$ is non-degenerate   on $\ker \alpha$, and $\coker(\mathcal{K} \hookrightarrow \mathbb{T}_A)\simeq \mathbb{T}_A/\mathcal{K} \simeq L[n]$ for a line bundle $L$ over $ \spec A$.
%\end{definition}

Suppose ${\bf X} = \spec A$, \ ${\bf Y} = \spec B$, with $A,B$ standard form cdgas. Recall that, in this case, we  have the identifications $\mathbb{L}_{{\bf X}}\simeq \mathbb{L}_A\simeq \Omega_A^1$ and $\mathbb{L}_{{\bf Y}}\simeq\mathbb{L}_B\simeq \Omega_B^1,$ along with Proposition \ref{defn_contact strc on good affines}. 

Let $(\rho ; \varLambda_{B} )$ be an isotropic structure  on $f: {\bf Y}\rightarrow ({\bf X}; \mathcal{K}\simeq \ker \alpha_A, L[n]) $, where  $f\simeq \spec \tilde{f}$ with the morphism  $\tilde{f}:A \rightarrow B$ of standard form cdgas;  \ $ \alpha_{A}: \mathbb{T}_A \rightarrow \mathcal{O}_A[n] $ is an $n$-contact form on $\spec A$, and the morphism $\rho$ is of the form $$ \rho: \mathbb{T}_B\hookrightarrow f^*(\ker \alpha_{A})\simeq B \otimes \ker \alpha_{A}, \text{ where} \quad \mathcal{O}(\spec B)=B.$$  
\begin{observation} \label{observ_some maps for big diagram}
%	We list some useful observations.
	\begin{enumerate}
	\item Using Proposition \ref{prop_pullback of cocone=cocone of the pullback} (with suitable modifications, if necessary), we have the weak equivalences $ f^*(\ker \alpha_A)\simeq \ker (f^*\alpha_A) $ and {$ f^*(\ker^{\vee} \alpha_A)\simeq \ker^{\vee} (f^*\alpha_A) $ },  hence the induced map $ \mathbb{T}_B\hookrightarrow \ker (f^* \alpha_A) $ again denoted by $\rho,$  with its dual map $\ker^{\vee} (f^*\alpha_A) \rightarrow \mathbb{L}_B. $ 
	\item If, in addition, the isotropic structure $\rho$ above is non-degenerate, then combining the equivalences in the previous item with Observation \ref{observation_non-deg isot str induces weak eqv on pullback}, we also get the weak equivalence $$\ker(f^*\alpha_A) \simeq  \ker^{\vee}(f^*\alpha_A) [n].$$
	%\item From the fiber sequence $\mathbb{T}_B\rightarrow f^*(\mathbb{T}_A) \rightarrow\mathbb{T}_f[1]$, there is a natural map $f^* (\ker \alpha_A) \rightarrow \mathbb{T}_f[1] $ as the first map in the sequence factors through $ f^* (\ker \alpha_A)$ via the isotropic structure $\rho$.
	
\end{enumerate}
\end{observation}

%\begin{observation}
%	Assume that $({\bf X}, \alpha_{\bf X})=({pt}, 0)$, where we consider  $\alpha_{\bf X}\equiv 0$ as a globally defined $k$-shifted contact form. Then any Legendrian structure $\varLambda_{{\bf Y}}$ on $f: {\bf Y}\rightarrow ({pt}, 0)$ gives rise to a $(k-1)$-shifted contact derived scheme $ ({\bf Y}, \varLambda_{{\bf Y}}) $, where the map $ \chi $ reduces to the contraction morphism $d_{dR}\varLambda_{{\bf Y}}\cdot : \mathbb{T}_{{\bf Y}} \rightarrow \mathbb{L}_{{\bf Y}}[k-1]$, which is a quasi-isomorphism on $\ker \varLambda_{{\bf Y}}$, and hence satisfies the contact condition in the sense of Definition \ref{defn_shiftedcontact}.
%\end{observation}

%\paragraph{Explicit affine models.}
\begin{construction} \label{Leg model}Let $k<0$ be an \bfem{odd integer such that $k\equiv 1\mod 4$}. It should be noted that the same setup will also work for more general case $k\not\equiv 3\mod 4$ by introducing extra $ (-1)^{(-i+1)(k+1)} $ factors in $d$ and $\phi$ below. 	This construction is a variation on  \cite[Example 3.3]{JS}.
	
	Our aim is to define an \bfem{affine   Legendrian structure} on a morphism $f$ of affine derived schemes, with shifted contact target. To this end, we construct standard form cdgas $A,B$ with a morphism $A\rightarrow B$ inducing $f$ and satisfying the desired properties.  
	%Consider % the case $k\in\mathbb{Z}_{<0} \text{ with } k\not\equiv 3\mod 4$. 
	%Assume $A, \alpha_0$ are in $k$-shifted contact Darboux form in the sense of Theorem \ref{contact darboux}. %For simplicity, we may further suppose  
	%This example is a variation 
\paragraph{Step-1: Setting a contact target $(\spec A, \alpha_0)$.} We essentially follow \cite[Example 2.21]{Berktav2024} in this step. Fixing $\ell=-\lfloor(k+1)/2\rfloor$\footnote{$ \ell=-k/2 $ if $k$ is even; otherwise $\ell= -(k+1)/2$. Here, $\lfloor - \rfloor$ denotes the \emph{floor function}.}, we start with a smooth $\K$-algebra $A^0:=A(0)$ of dimension $m_0$, having degree 0 variables $x_1^0,  \dots, x_{m_0}^0$ in $A(0)$ such that $d_{dR}x_1^0,\dots,d_{dR}x_{m_0}^0 $ form a  basis for $\Omega_{A(0)}^1$. %The rest follows \cite[Example 2.21]{Berktav2024}. %This choice can be made by localizing $A(0)$ if necessary.

Next, choosing non-negative integers $m_1,\dots, m_{\ell}$, we define the  variables $x_j^{-i}, y_j^{k+i}, z^k \in A$ so that $A$ is a \emph{commutative graded algebra} freely generated over $A(0)$ by these variables, and $\Omega^1_{A}$ is  the free  $A$-module of finite rank with basis $ \{d_{dR}x_j^{-i}, d_{dR}y_j^{k+i}, d_{dR}z^k : i= 0, 1, \dots, \ell, \ j= 1,2, \dots, m_i\}.$
	
	To make $A$ a \emph{cdga}, we define the \bfem{internal differential $d$ on $A$} as \begin{align} \label{defn_internal d contact}
			d|_{A(0)}&=0; \ dx_j^{-i} =  \dfrac{\partial H}{\partial y_j^{k+i}} \text{ for all } i>0,j; \  \ dy_j^{k+i} =  \dfrac{\partial H}{\partial x_j^{-i}} \text{ for all } i,j; \text{ and } \nonumber \\ -kdz^k&= H+d\Big[\sum_{i,j} (-1)^{i} i x_j^{-i}y_j^{k+i} \Big],
		\end{align} where $H$ is the \emph{Hamiltonian} satisfying the ``classical master equation". \cite[Example 2.21]{Berktav2024} 
	%We then let  $\alpha_0 \in (\Omega_A^1)^k $  such that 
%	$ d_{dR}\alpha_0= \sum_{i,j}^{ } d_{dR}x_j^{-i} d_{dR}y_j^{k+i}$   satisfying
%	$dH=0, \ d_{dR}H+d\phi=0, \text{ and } d_{dR}\phi=kd_{dR}\alpha_0$; and 
	%Moreover, from Theorem \ref{contact darboux}, we may assume w.l.o.g. that 
	    also shows that the element \[\alpha_0 = d_{dR}z^k + \sum_{i=0}^{\ell} \sum_{j=1}^{m_i} y_j^{k+i}\dR x_j^{-i} \in \Omega_A^1 [k]\] defines a \bfem{$k$-contact form}  such that  one has the natural splitting 
	\begin{equation} 
	\mathbb{T}_A|_{\spec H^0(A)}=  \ker \alpha_0 |_{\spec H^0(A)} \oplus Rest|_{\spec H^0(A)}
	\end{equation} where 
	$ Rest|_{\cspec H^0(A)} =\big\langle \partial/\partial z^k \big\rangle_{A(0)}.$
	 Moreover, it has been shown in \cite{kib,Berktav2024} that the pair \[\Big(H,\phi:=k\alpha_0-kd_{dR}z^k-\dR\Big[\sum_{i,j} (-1)^{i} i x_j^{-i}y_j^{k+i} \Big]\Big),\] with $H,\alpha_0$ as above, provides a solution to the  equations \begin{equation} \label{important relations}
		dH=0 \text{ in } A^{k+2}, \ d_{dR}H+d\phi=0 \text{ in } (\Omega^1_{A})^{k+1}, \text{ and } d_{dR}\phi=kd_{dR}\alpha_0,
	\end{equation} which implies  $d\alpha=0$, such that the element $\phi$ above can  be written explicitly as $$\phi:=k\alpha_0-kd_{dR}z^k-\dR[\cdots]=  \sum_{i=0}^{\ell} \sum_{j=1}^{m_i} \big[(-i)x_j^{-i}d_{dR}y_j^{k+i}+(k+i)y_j^{k+i}d_{dR}x_j^{-i}\big].$$

Note that using the terminology from \cite{JS}, we can decompose\footnote{Due to the degree reasons, $H$ must be at most linear in $ y_j^{k+i}. $} $H$ as 
\begin{equation}
	H=H_+ +  \sum_{i,j} H^{(-i+1)}_jy_j^{k+i}, \text{ with } H_+\in A_+^{k+1}, \ H^{(-i+1)}_j\in A_+^{-i+1}.
\end{equation}Likewise, setting the element\footnote{Actually, $\phi_+:= -\sum_{i=0}^{\ell} \sum_{j=1}^{m_i} (-1)^{(-i+1)(k+1)} y_j^{k+i}d_{dR}x_j^{-i}$, but the term $(-1)^{(-i+1)(k+1)}=1$ when $k$ is odd.}$\phi_+:=\dR z^k-\alpha_0= -\sum_{i=0}^{\ell} \sum_{j=1}^{m_i}  y_j^{k+i}d_{dR}x_j^{-i}$, Joyce and Safronov showed in \cite[Remark 2.15 \& Example 3.3]{JS} that the pair $(H_+, \phi_+)$ satisfies the equations \begin{equation}\label{important relations+}
dH_+=0, \ d_{dR}H_+ + d\phi_+ =0, \text{ and } d_{dR}\phi_+=-\dR\alpha_0.
\end{equation}%For details, we refer to \cite[Remark 2.15 \& Example 3.3]{JS}. Also, see Appendix \ref{section_models for LAGNBH thm}.

Before the next step, it will be useful to express $\alpha_0$ and $dz^k$ in terms of the pair $(H_+, \phi_+)$ satisfying  (\ref{important relations+}). By definition, we can simply write
\begin{equation}\label{alph_0 and dz in term of H_+ phi_+}
	\alpha_0= \dR z^k - \phi_+, \quad \text{and} \quad dz^k=H_+.
\end{equation} The second equality in (\ref{alph_0 and dz in term of H_+ phi_+}) follows from combining $d\alpha_0=0$ with (\ref{important relations+}). In fact, one computes $0=d\alpha_0= d (\dR z^k - \phi_+) = -\dR \circ dz^k-d\phi_+= -\dR \circ dz^k + \dR H_+$, leading to $dz^k=H_+.$
%With the above setup in hand, we now give an \bfem{affine model for Legendrian structures} on morphisms $f$ into $\spec A$. To this end, we construct another standard form cdga $B$ with a morphism $A\rightarrow B$ inducing $f$ and satisfying the desired properties.  
% \[s = 
% \begin{cases*} -\ell, & $ k $ \mbox{ even}, \\
% -\ell-1, & $ k $ \mbox{ odd}.
% \end{cases*}
% \]
\paragraph{Step-2: Setting the source $\spec B$.} Let  $B^0:=B(0)$ be a smooth algebra of dimension $m_0+n_0$ and $\beta^0:A^0\rightarrow B^0$ a smooth morphism. Assume that there exist  $u^0_1, \dots, u^0_{n_0}$ in $B^0$ such that $\{d_{dR}\tx_{1}^0, \dots, \tx^0_{m_0}, d_{dR}u^0_1, \dots, d_{dR}u^0_{n_0} \}$ is a basis over $B^0$ for $\Omega_{B^0}^1$, with $\tx^0_j=\beta^0(x^0_j)$. 

 Let $\beta_+:A_+ \rightarrow B$ be a morphism of \emph{commutative graded algebras} given by $\beta_+|_{A^0}=\beta^0$ and $\tx^{-i}_j=\beta_+(x^{-i}_j)$. Writing $s:=-\lfloor k/2\rfloor$\footnote{$s=\ell$ if $k$ is even; otherwise $s=\ell+1$.}, we define the \emph{commutative graded algebra} $B$ to be the free graded algebra over $B(0)$ generated by the variables  
\begin{align} 
& \tx_1^{-i}, \tx_2^{-i}, \dots, \tx_{m_i}^{-i} & & \text{ in degree } (-i) \text{ for } i= 1, 2, \dots, \ell, \nonumber \\
&u_1^{-i}, u_2^{-i}, \dots, u_{n_i}^{-i} & & \text{ in degree } (-i) \text{ for } i=1, 2, \dots, s, \nonumber \\
&v_1^{k-1+i}, v_2^{k-1+i}, \dots, v_{n_i}^{k-1+i} & & \text{ in degree } (k-1+i) \text{ for } i=0,\dots, s.
%&v_1^{k}, v_2^{k}, \dots, v_{n_1}^{k} & & \text{ in degree } k.\nonumber \\
%&v_1^{k-1}, v_2^{k-1}, \dots, v_{n_0}^{k-1} & & \text{ in degree } (k-1).
\end{align}%where $  w^{k-1} $ is the \bfem{distinguished generator} in degree $k-1$. 

Now, choose an element $G \in B^k$, the \emph{superpotential},  satisfying the \textit{relative version} of  classical master equation (\ref{relative CME}), which determines the \bfem{differential $d_B$ on $B$} as (cf. Equation (\ref{differantial on B})): \begin{align*}
	&d_B|_{B^0}=0, \\ &d_B\tx_j^{-i}=(-1)^{1-i} \beta_+(H_j^{-i+1}), \quad \text{ for } i= 1, 2, \dots, \ell, \\ &d_Bu_j^{-i}=(-1)^{(1-i)k} \partial G/ \partial v_j^{k-1+i}, \quad d_Bv_j^{k-1+i}= \partial G/ \partial u_j^{-i} \quad \text{ for } i=0, 1, 2, \dots, s.
\end{align*} \cite[Example 3.3]{JS} shows that  $d^2_B=0$ on each generator above, implying that $(B, d_B)$ is a standard form cdga.

\paragraph{Step-3: Setting  the map into $\spec A$.} We now extend $\beta_+$ to a morphism $\beta:A\rightarrow B$ of \emph{commutative graded algebras} via $\beta|_{A_+}=\beta_+$ and setting
\begin{align*}
	\beta(y_j^{k+i})&= (-1)^{-i+1} \partial G/\partial \tx_j^{-i},\\
(k-1)\beta(z^k)&=G+\sum_{i,j} i \tx^{-i}_j \partial G/\partial \tx_j^{-i} \\ &+ d_B\Big[\sum_{i,j} (-1)^{(1-k-i)} (1-k-i) (-1)^{(1-i)k}v_j^{k+i-1}u_j^{-i}\Big].
\end{align*}%We then get the map $\spec \beta: \spec B \rightarrow (\spec A,\alpha_0)$, with the  differential on $B$ above. %defined by
%\begin{align}
	%-kd_B\tilde{z}^k&=\beta(H)+d_B\Big[\sum_{i,j} (-1)^{i} i\tx_j^{-i} \partial G/\partial \tx_j^{-i} \Big]\\
%-(k-1)d_Bw^k&=G + (1-k)\beta(z^k)+ \sum_{i,j}  i\tx_j^{-i} \partial G/\partial \tx_j^{-i}\\ \nonumber &+ d_B\Big[\sum_{i,j} (-1)^{(1-k-i)} (1-k-i) (-1)^{(1-i)k}v_j^{k+i-1}u_j^{-i}.\Big] 
%\end{align} 
To promote $\beta:A\rightarrow B$ a morphism of \emph{cdgas}, it remains to show that  $d_B \circ \beta = \beta \circ d$ for the variables $x_j^{-i}, y_j^{k+i}, z^k \in A$.   Let us clarify:

\begin{itemize}
	\item From \cite[Example 3.3]{JS}, the equality $d_B \circ \beta = \beta \circ d$ holds true for each $ x_j^{-i}, y_j^{k+i}.$
	\item Using the previous item and the definition of $\beta$ on $ x_j^{-i}, y_j^{k+i}$, we get \begin{align*}
		(k-1)d_B \circ \beta (z^k) &= d_BG + d_B \Big[\sum_{i,j} i \tx^{-i}_j \partial G/\partial \tx_j^{-i}\Big] +d_B^2[\cdots] \\
		&= -\beta(H+H_+)  + \sum_{i,j} i \beta(dx^{-i}_j) (-1)^{1-i}\beta(y^{k+i}_j) \\ &+ \sum_{i,j} i \beta(x^{-i}_j) (-1)^{-i}(-1)^{1-i}\beta(dy^{k+i}_j)+0\\
		&= -\beta(H+H_+)  - \beta\Big(\sum_{i,j} i dx^{-i}_j (-1)^{i} y^{k+i}_j\Big) - \beta\Big(\sum_{i,j} i x^{-i}_j dy^{k+i}_j\Big)\\
		&= -\beta(H)- \beta\circ d\Big( \sum_{i,j} (-1)^{i} i x_j^{-i}y_j^{k+i} \Big)-\beta(H_+)\\
		&= -\beta\Big[\ \underbrace{H+ d\Big( \sum_{i,j} (-1)^{i} i x_j^{-i}y_j^{k+i}\Big)}_{=-kdz^k} \ \Big] -\underbrace{\beta(H_+)}_{=\beta(dz^k)} \\
		&= (k-1) \beta \circ d (z^k).
	\end{align*}
\end{itemize}We  then conclude that $d_B \circ \beta (z^k)= \beta \circ d(z^k)$ as well; hence, both $\beta, \beta_+$ are \emph{morphisms of cdgas} such that $\beta_+$ is a submersion. Thus, we get the map $\spec \beta: \spec B \rightarrow (\spec A,\alpha_0)$, as desired.

\paragraph{Step-4: Isotropic structure.} We next define $\varLambda \in (\Omega_B^1)^{k-1}$ to be the element of the form
\begin{equation} \label{path btw 0 and f^*alpha}
	\varLambda = \sum_{i=0}^{s} \sum_{j=1}^{n_i}  u_j^{-i}\dR v_j^{k-1+i}.
\end{equation}Joyce and Safronov \cite{JS} show that the element $ d_{dR}\varLambda= \sum_{i=0}^{s} \sum_{j=1}^{n_i} d_{dR}u_j^{-i}d_{dR}v_j^{k-1+i} $ 
 %a non-degenerate isotropic structure $h:=(h^0, 0, \dots)$ on $\spec \beta: \spec B \rightarrow (\spec A, \omega_{can})$, with the equations $dh^0= \beta_* (\omega^0), \  d_{dR}h^0=0,$
 satisfies the equations 
\begin{equation}\label{defining relations for B,d_B}
d_BG=-\beta(H+H_+), \ d_{dR}G+d_B\psi=-\beta_*(\phi+\phi_+), \text{ and } d_{dR}\psi=(k-1)d_{dR}\varLambda,
\end{equation} where  $ d_{dR}\varLambda, d_{dR}\alpha_0 $ play the roles of $ h^0, \omega^0$ in Construction (\ref{Lag model}), respectively;  and $ H_+, \phi_+, \psi $ are as in Construction (\ref{Lag model}). %In addition to that, \cite{JS} also proves that $d_B(\dR\varLambda)=\beta_* (d_{dR}\alpha_0)$, 
Thus, we have the following observations:

\begin{itemize}
	\item As $\psi:= \displaystyle \sum_{i=0}^{s} \sum_{j=1}^{n_i} \big[-iu_j^{-i}d_{dR}v_j^{k-1+i} +(-1)^{(1-i)k}(k-1+i)v_j^{k-1+i}d_{dR}u_j^{-i}\big]$, one can write \[(k-1)\varLambda=\psi + \dR \Big[\sum_{i,j} (-1)^{(1-k-i)} (1-k-i) (-1)^{(1-i)k}v_j^{k+i-1}u_j^{-i}\Big]. \]
	\item Note that $ \phi + \phi_+ = (k-1)\alpha_0 + (1-k) \dR z^k - \dR\Big[\sum_{i,j} (-1)^{i} i x_j^{-i}y_j^{k+i} \Big].$
	\item Using the definitions of $d_B$ and $\beta$ together with the equations (\ref{defining relations for B,d_B}), one computes  
	\begin{align*}
		(k-1) d_B\varLambda&= d_B\psi + d_B\circ \dR [\cdots] \\
		&= -\dR G - \beta_*(\phi + \phi_+) + d_B\circ \dR [\cdots]\\
		&= -\dR G +(1-k) \beta_*(\alpha_0) \ + (k-1) \dR \beta (z^k) \\ &- \dR\Big[\sum_{i,j}  i \tx_j^{-i}\partial G/\partial \tx_j^{-i} \Big] + {d_B\circ \dR [\cdots]}\\
		&= (1-k)\beta_*(\alpha_0),
	\end{align*} where $(k-1) \dR \beta (z^k)= \dR G+\dR \sum_{i,j} i \tx^{-i}_j \partial G/\partial \tx_j^{-i}+{\dR \circ d_B [\cdots]} $ is used in the third line to get cancellations. It then follows that 	
	\[ d_B \varLambda = -\beta_* (\alpha_0), \text{ implying also that }  d_B(\dR\varLambda)=\beta_* (d_{dR}\alpha_0).\] %See also \cite[Example 3.3]{JS}.
\end{itemize}
%; and such that %; and $\psi$ is given by  Equation (\ref{defn_psi}). 
%As before,  we may also assume w.l.o.g. that 
%\begin{align*}
%\ker \varLambda |_{\cspec H^0(B)} &=\big\langle \partial/\partial u_1^{-i}, \dots, \partial/\partial u_{n_i}^{-i},  \partial/\partial v^{k-1+i}_1, \dots, \partial/\partial v^{k-1+i}_{n_i} : i=0,1, \dots, s \ \big\rangle_{B(0)},  \\
%Rest|_{\cspec H^0(B)} &=\big\langle \partial/\partial U_1^{0}, \partial/\partial\tilde{X}_1^0, \partial/\partial\tx_1^{-i},  \dots, \partial/\partial\tx_{m_i}^{-i} :  i= 1, 2, \dots, \ell\big\rangle_{B(0)} . 
%\end{align*} 
From the computation above, we conclude that $ d_{dR}\varLambda $  satisfies the  defining equations, see (\ref{simplification for h}), $d_B(d_{dR}\varLambda)= \beta_* (d_{dR}\alpha_0), \  d_{dR}(d_{dR}\varLambda)=0$; and hence, $(d_{dR}\varLambda, 0, \dots)$, with $\varLambda \in (\Omega^1_B)^{k-1}$, defines an \emph{isotropic structure} from $0$ to $ f^*(d_{dR}\alpha_0) $. Likewise, since $-\varLambda \in (\Omega^1_B)^{k-1}$ is such that $d_B (-\varLambda) = \beta_* (\alpha_0)$, the element $ -\varLambda $ defines a \emph{path between $0$ and $f^*(\alpha_0)$. }

Note the  path $ -\varLambda: 0 \rightsquigarrow f^*(\alpha_0)$ is part of contact isotropic data. It still remains to find a suitable morphism $\rho: \mathbb{T}_B \rightarrow f^*(\ker \alpha_0)$ commuting Diagram \ref{defn_isotropic for leg diagram}. To set $\rho: \mathbb{T}_B \rightarrow f^*(\ker \alpha_0)$, it is enough to consider the map  $$\rho\otimes_{B}H^0(B):\mathbb{T}_B\otimes_{B}H^0(B)\rightarrow \ker \alpha_0\otimes_{A}H^0(B),$$ which sends\footnote{In fact, we use the map $\langle \partial/\partial \tx^{-i}_j \rangle_{H^0(B)} \mapsto \langle \partial/\partial x^{-i}_j - y_j^{k+i}\partial/\partial z \rangle_{H^0(B)} \subset f^*(\ker \alpha_0)$ combined with the identification(s) via the quotient map $\ker \alpha_0 \simeq \mathbb{T}_A/\langle \partial/\partial z\rangle_{A(0)} \simeq \langle \partial/\partial x^{-i}_j, \partial/\partial y_j^{k+i} \rangle_{A(0)}$.} $\langle \partial/\partial \tx^{-i}_j \rangle_{H^0(B)} \mapsto \langle \partial/\partial x^{-i}_j \rangle_{H^0(B)}$, with trivial action on the other vector fields. By construction, such morphism $\rho$ makes
Diagram \ref{defn_isotropic for leg diagram}
commute; hence, the \bfem{pair $(\rho, \varLambda)$}, where \begin{itemize}
	\item $ -\varLambda: 0 \rightsquigarrow f^*(\alpha_0)$ is the path defined by the element $\varLambda \in (\Omega^1_B)^{k-1}$ in (\ref{path btw 0 and f^*alpha}),
	\item  $\rho: \mathbb{T}_B \rightarrow f^*(\ker \alpha_0) \text{ is the morphism, sending } \langle \partial/\partial \tx^{-i}_j \rangle \mapsto \langle \partial/\partial x^{-i}_j \rangle$,	
\end{itemize}defines a \bfem{contact isotropic structure} on the map  $f:=\spec \beta: \spec B \rightarrow (\spec A;\alpha_0, L[k]).$

\paragraph{Step-5: Non-degeneracy.} %??? From definitions, $ d_{dR}\varLambda $ makes Diagram (\ref{isotropic strc_Leg}) commute (we simply use $ \varLambda $ instead of $ \varLambda_B $), and hence it remains to show that $ d_{dR}\varLambda $ is non-degenerate. To this end, we need to show that the morphism $\chi|_{\ker \varLambda}: 	\mathbb{T}_{f}|_{\ker \varLambda} \rightarrow \mathbb{L}_{B}|_{(\ker \varLambda)^{\vee}}[k-1]$ at the bottom face of the diagram (\ref{isotropic strc_Leg}) is a quasi-isomorphism, where $f:=\spec \beta: \spec B \rightarrow (\spec A,\alpha_0)$ carrying the isotropic structure $(d_{dR}\varLambda, 0, \dots)$ as above. In fact, it suffices to apply $-\otimes_{B}H^0(B)$ to Diagram (\ref{isotropic strc_Leg}) and show that the corresponding map $\chi|_{\ker \varLambda} \otimes_{B}H^0(B)$ is an isomorphism of complexes of $H^0(B)$-modules. 

So far, we have obtained in the last step a contact isotropic structure $(\rho, \varLambda)$ on the map  $f:=\spec \beta: \spec B \rightarrow (\spec A;\alpha_0, L[k]).$ It remains to show that this contact isotropic structure is non-degenerate. 
To this end, we need the induced morphism $ \mathbb{T}_{\rho} \rightarrow \mathbb{L}_{B}[k-1]$ to be a quasi-isomorphism, where $\mathbb{T}_{\rho}:=\mathrm{hofib}\left(\rho: \mathbb{T}_{B}\rightarrow f^*(\ker \alpha_0)\right)$. As for the case of Lagrangians \cite{JS}, it suffices to apply $-\otimes_{B}H^0(B)$ and show the morphism  \[\mathbb{T}_{\rho}\otimes_{B}H^0(B) \rightarrow \mathbb{L}_{B}[k-1]\otimes_{B}H^0(B)\] is an isomorphism of complexes of $ H^0(B)$-modules.

Note that %we are working with sufficiently nice cdgas; and hence, in what follows, we use $\Omega^1_{(-)}$ as a model for $\mathbb{L}_{(-)}$. 
since $\mathbb{T}_B\xrightarrow{\pi} f^*\mathbb{T}_A \rightarrow \mathbb{T}_f[1]$ is an exact sequence with the cofiber $\mathbb{T}_f[1]$, we can use the \emph{cone} of the morphism $\pi=(\Omega_{\beta}^1)^{\vee}: (\Omega_B^1)^{\vee}\rightarrow (\Omega_A^1)^{\vee}$, induced by the map $\beta:(A,d)\rightarrow (B,d_B),$ as a model for $ \mathbb{T}_{f} \otimes_{B}H^0(B)$. Recall that for such cdgas, $Cone((\Omega_{\beta}^1)^{\vee})$ can be given as %$ \mathbb{T}_{B/A} \otimes_{B}H^0(B)$
\begin{equation}
Cone((\Omega_{\beta}^1)^{\vee})=\Bigg( \big((\Omega_B^1)^{\vee} \otimes_{B}H^0(B)\big) \oplus   \big((\Omega_A^1)^{\vee} [-1]\otimes_{A}H^0(B)\big) ,       \begin{pmatrix}
d_B & 0 \\
(\Omega_{\beta}^1)^{\vee} & d
\end{pmatrix}   \Bigg).
\end{equation} Likewise, for the fiber-cofiber sequence $ \mathbb{T}_{B}\xrightarrow{\rho} f^*(\ker \alpha_0) \rightarrow \mathbb{T}_{\rho}[1]$, we use the following complex as a model for $ \mathbb{T}_{\rho} \otimes_{B}H^0(B)$:
\begin{equation}
\Bigg( \big(\mathbb{T}_{B} \otimes_{B}H^0(B)\big) \oplus   \big(\ker \alpha_0 [-1]\otimes_{A}H^0(B)\big) ,       \begin{pmatrix}
		d_B & 0 \\
		\rho & d
	\end{pmatrix}   \Bigg).
\end{equation}

Let us write $ \mathbb{T}_{\rho} \otimes_{B}H^0(B)$ and $ \Omega_B^1[k-1] \otimes_{B}H^0(B) $ as direct sums of $H^0(B)$-modules. For each $i\geq 0$, we have:
\begin{align*}
\left(\mathbb{T}_{\rho} \otimes_{B}H^0(B)\right)^i&= \langle \partial/\partial \tx^{-i}_j\rangle \oplus \langle \partial/\partial u_j^{-i}, \partial/\partial v_{j}^{-i} \rangle \oplus
\ \langle \partial/\partial x^{-i+1}_j\rangle \oplus \langle \partial/\partial y^{1-i}_j\rangle, \nonumber \\
\left(\Omega_B^1[k-1] \otimes_{B}H^0(B)\right)^{-i}&= \langle d_{dR}\tx^{k-i+1}_j\rangle\oplus \langle d_{dR}u_j^{k-1+i}, d_{dR} v_j^{k-1+i}\rangle.
\end{align*}From \cite[Example 3.3]{JS}, the  map $\chi_{\rho}: 	\mathbb{T}_{\rho} \otimes_{B}H^0(B)\rightarrow \Omega_B^1[k-1] \otimes_{B}H^0(B)$ is then given degree-wise by 
\begin{equation}
 \chi_{\rho}^{*}=\begin{pmatrix}
0 & 0 & \star & \beta_* \circ d_{dR}\alpha_0 |_{\ker \alpha} \ \cdot \\
0 &  d_{dR}\varLambda \ \cdot & \star & 0
\end{pmatrix} ,
\end{equation} where $ \star$'s are some arbitrary morphisms; $\beta_* \circ d_{dR}\alpha_0 \cdot \text{ maps } \partial/\partial y^{\mu}_j \mapsto d_{dR}x^{k-\mu}_j \mapsto d_{dR}\tx^{k-\mu}_j$; and $d_{dR}\varLambda \ \cdot$ maps $\langle \partial/\partial u_j^{\mu}, \partial/\partial v_j^{\mu} \rangle \mapsto \langle d_{dR}u_j^{k-1-\mu},d_{dR} v_j^{k-1-\mu} \rangle$, all up to sign. Clearly, both such maps in $\chi_{\rho}^{\mu}$, for each $\mu$, are isomorphism of $H^0(B)$-modules.

\cite[Example 3.3]{JS} also shows that there is a  subcomplex $(C,d_B)$ of $ \mathbb{T}_{\rho} \otimes_{B}H^0(B) $, with $C=\{0\} \oplus \langle \partial/\partial u_j^{-i}, \partial/\partial v_{j}^{-i} \rangle \oplus \{0\} \oplus \langle \partial/\partial y^{1-i}_j\rangle$, such that the inclusion $\iota: C \hookrightarrow \mathbb{T}_{\rho} \otimes_{B}H^0(B)$ is a quasi-isomorphism and that $\chi_{\rho}|_{C}:=\chi_{\rho}\circ\iota: C \rightarrow \Omega_B^1[k-1] \otimes_{B}H^0(B)$ is a strict isomorphism of complexes; and hence, $\chi_{\rho}$ is a quasi-isomorphism.

%Thus, $\chi $ maps $ \ker \varLambda |_{\cspec H^0(B)} $ isomorphically onto its image, where we consider $ \ker \varLambda |_{\cspec H^0(B)} $ as a subcomplex of $ \mathbb{T}_{f} \otimes_{B}H^0(B) $ of the form $\{0\} \oplus \langle \partial/\partial u_j^{-i}, \partial/\partial v_{j'}^{-i'} \rangle \oplus \{0\} \oplus \{0\}$. 

As a result,  the contact isotropic structure $ (\varLambda, \rho) $ on $\spec \beta$ is non-degenerate; and hence, $\spec B$ is Legendrian in $(\spec A; \alpha_0, L[k]).$

\end{construction}

We complete the prototype construction of an affine model for Legendrian structures, leading to the following definition.

\begin{definition} \label{defn_Leg model}
	Let $(\spec B, \varLambda, \rho)$ be a Legendrian in $(\spec A, \alpha_0)$ as in Construction \ref{Leg model}. We then say that the data $(A, \alpha_0, B, \beta, \varLambda,\rho)$ is in  \bfem{Legenderian-Darboux form}.
\end{definition}	
\subsection{A Legendrian-Darboux theorem for contact derived schemes}\label{section_LEG-DA-tthm}	In brief, we would like to show that Legedrians $f: \bf Y \rightarrow {\bf X}$ in negatively shifted contact derived schemes are locally modeled on standard Legendrian-Darboux forms $(A, \alpha_0, B, \beta, \varLambda,\rho)$ as in Construction \ref{Leg model}. More precisely, we have:
\begin{theorem} \label{Legendrian-Darboux}
		Let $({\bf X}; \mathcal{K}, \kappa, L[k])$ be a (locally finitely presented) $k$-shifted contact derived $\mathbb{K}$-scheme, with $k\in\mathbb{Z}_{<0}$, and $f: \bf Y \rightarrow {\bf X}$ be a morphism of derived $\mathbb{K}$-schemes carrying a Legendrian structure $ (\varLambda', \rho').$ Given $y\in \bf Y$, $x\in \bf X$ with $f(y)=x$, let $(A, i: \spec A \hookrightarrow \bf X),$ with $p\in \cspec H^0(A)$, be a standard form affine neighborhood of $x$ such that $i(p)=x$ (as in Theorem \ref{localmodelthm}). 
		Then  there exist \begin{enumerate}
			\item [(1)] a standard form cdga $B$ and a point $q\in \spec B,$ \item [(2)] a morphism $\beta: A \rightarrow B$ of cdgas with $\spec \beta: q\mapsto p$, such that  $\beta_+:=\beta|_{A_+}: A_+ \rightarrow B$ is a submersion minimal at $q$, \item [(3)] a Zariski open inclusion $j:\spec B \hookrightarrow \bf Y$, with $j(q)=y$, such that the  diagram 
		\begin{equation} 
		\begin{tikzcd}
			\spec B \arrow[d, "\spec \beta"'] \arrow[r, "j", hook] & \bf Y \arrow[d, "f"] \\
			\spec A \arrow[r, "i"', hook]                          & \bf X               
		\end{tikzcd}
		\end{equation} commutes up to homotopy, and \item [(4)] the pullback of the Legendrian structure  to $ \spec \beta $ satisfying $j^*(\varLambda')\simeq \varLambda$ and $\rho=\rho'|_{\mathbb{T}_B}$ such that $(\spec B, \varLambda,\rho)$ is a Legendrian in $(\spec A, \alpha_0)$, where $(A, \alpha_0, B, \beta, \varLambda, \rho)$ is in Legenderian-Darboux form in the sense of Definition \ref{defn_Leg model}. 
		\end{enumerate}
\end{theorem} 

\emph{Disclaimer.}  In what follows, we will consider the case of \bfem{odd negative shifts} only, %$``k<0 \text{ with } k\not\equiv 3\mod 4"$ \emph{only}, 
for which Construction \ref{Leg model} gives the desired explicit model. Our result is a variation on \cite[Theorem 3.7]{JS}. %Therefore, There will be possible extra simplifications throughout the proof. Modifications for 
The other cases can be dealt with by adopting almost the same modifications as in \cite[Example 3.5; Section 4]{JS}; hence, the main body of the proof will remain the same. The modifications are relatively straightforward, but sometimes cumbersome. We will leave details to the reader.  

\pf Let $y\in \bf Y$, $x\in \bf X$ with  $f(y)=x$. For $x\in \bf X$, apply Theorem \ref{contact darboux} to find an affine contact neighborhood    $\big(A,  \alpha_0 \big)$ of $x$, with $p \in \spec H^0(A)$, where $i: \spec A \hookrightarrow \bf X$ is an open inclusion mapping $i: p\mapsto x$; \ $A$ is a standard form that is minimal at  $p$; and $\alpha_0$ is a standard $k$-shifted contact form on $\spec A$ such that $A, \alpha_0$ are in Darboux form.

\paragraph{(0) Representing the Darboux form $A, \alpha_0$ in local variables.} Let $\ell=-\lfloor(k+1)/2\rfloor$, then Theorem \ref{contact darboux} implies that we can choose a smooth $\K$-algebra $A(0)$ of dimension $m_0$ with  generators $x_1^0, x_2^0, \dots, x_{m_0}^0$ of degree 0, and a set of variables $x_j^{-i}, y_j^{k+i}, z^k \text{ in } A$  such that  $\Omega^1_{A}$ is  the free  $A$-module of finite rank with basis $ \{d_{dR}x_j^{-i}, d_{dR}y_j^{k+i}, d_{dR}z^k : i= 0, 1, \dots, \ell, \ j= 1,2, \dots, m_i\}.$ 

Here, due to the construction, except the distinguished generator $z^k$ of degree $k$, all other generators are defined as symplectic pairs via the non-degenerate pairing $\ker \alpha_0  \rightarrow \ker^{\vee} \alpha_0 [k]$. Thus, the quotient ${\mathbb{T}_A}/{\ker \alpha}$ can  be identified with a $k$-shifted line bundle $Rest$ (viewed as a perfect complex concentrated in degree $-k$) such that one has the  splitting 
 \begin{align*}
	\ker \alpha_0 |_{\cspec H^0(A)} &=\big\langle \partial/\partial x_1^{-i}, \dots, \partial/\partial x_{m_i}^{-i},  \partial/\partial y^{k+i}_1, \dots, \partial/\partial y^{k+i}_{m_i} : i=0,1, \dots, \ell \ \big\rangle_{A(0)},  \\
	Rest|_{\cspec H^0(A)} &=\big\langle \partial/\partial z^k \big\rangle_{A(0)} . 
\end{align*} 

Such splitting implies that we have	$ d_{dR}\alpha_0= \sum_{i,j}^{ } d_{dR}x_j^{-i} d_{dR}y_j^{k+i}$, along with the pair $(H,\phi)$ satisfying certain equations (cf. Construction \ref{Leg model}). Moreover, scaling $z^k$ we may assume $\iota_{\partial/\partial z^k} \alpha =1$. Then \cite[Theorem 3.13]{kib} proves that the corresponding conditions above uniquely (up to interchange of $x_{j}^{-i} \text{ and } y_{j}^{k+i}$) determine the representation of the element $\alpha_0\in \Omega_A^1[k]$ as 
\begin{equation*}
	\alpha_0 = d_{dR}z^k +  \sum_{i,j} y_j^{k+i}d_{dR}x_j^{-i}.
\end{equation*}Using  the  equations
 $dH=0, \ d_{dR}H+d\phi=0, \text{ and } d_{dR}\phi=kd_{dR}\alpha_0$, \cite{kib} also shows  that the differential $d$ is fully given as in  (\ref{defn_internal d contact}) and that $(A,d)$ is identified  with the cdga in Construction \ref{Leg model}.

\paragraph{(1) Constructing a local representative for $f$.} Let us start with  the pullback diagram 
\begin{equation} 
\begin{tikzcd}
	{\bf Y} \times_{\bf X} \spec A \arrow[dr, phantom, "\lrcorner", near start]           \arrow[d, "\pi_1"'] \arrow[r, "\pi_2"] & \spec A \arrow[d, "i", hook] \arrow[rr, "\spec(\iota)"] &  & \spec A_+ \\
	\bf Y \arrow[r, "f"]                                      & \bf X                                             &  &          
\end{tikzcd}
\end{equation} along with the map $\spec (\iota): \spec A \rightarrow \spec A_+$ induced by $\iota: A_+ \hookrightarrow A$, where $A_+$ (and other relevant objects like $H_+, \beta_+, \phi_+$) is/are as in Construction \ref{Lag model}.

Applying Theorem \ref{thm_nice models for morphisms} to the map $\spec (\iota) \circ \pi_2:{\bf Y} \times_{\bf X} \spec A \rightarrow \spec A_+$, with the map $id_{\spec A_+}  \text{ and the points } z \in {\bf Y} \times_{\bf X} \spec A,   p \in \spec A_+$, we get  a standard form cdga $B$, a point $q\in \spec B,$ a submersion $\beta_+: A_+ \rightarrow B$ minimal at $q$ with $\spec \beta_+ (q)=p$, and a Zariski open inclusion $j':\spec B \hookrightarrow {\bf Y} \times_{\bf X} \spec A$, with $j'(q)=y$, such that the following diagram  commutes:
\begin{equation} 
\begin{tikzcd}
	\spec B  \arrow[dd, "\spec \beta_+"', dashed] \arrow[r, "j'", dashed] \arrow[rrd, "\pi_2 \circ j'" description, dotted, near start] & {\bf Y} \times_{\bf X} \spec A \arrow[rd, "\pi_2"] \arrow[dd, "\spec (\iota) \circ \pi_2" description] &                                    \\
	&                                                                                                        & \spec A \arrow[ld, "\spec(\iota)"] \\
	\spec A_+ \arrow[r,"id"]                                                                                              & \spec A_+                                                                                              &                                   
\end{tikzcd}
\end{equation} Here, the morphism $\pi_2\circ j': \spec B \rightarrow \spec A$ is induced by a map $\beta^{\infty}: A\rightarrow B$ in $cdga_{\K}^{\infty}$, with $\beta^{\infty}\circ \iota\simeq \beta_+$; \ i.e., \ $\pi_2\circ j'\simeq \spec \beta^{\infty}.$ 

From \cite[$\S 4.2$ ]{JS}, we can conclude that since $A$ is free over $A_+$ (hence cofibrant over $A_+$), up to equivalence, $ \beta^{\infty} $ descends to a morphism $\beta:A\rightarrow B$ in $cdga_{\K}^{\leq 0}$ such that $\beta \circ \iota = \beta_+$; and thus, $\pi_2\circ j'\simeq \spec \beta.$ That is, we obtain the  diagram 

\begin{equation}
	\begin{tikzcd}
		& \spec B \arrow[rr, "j'"] \arrow[rd, "j:=\pi_1\circ j'" description] \arrow[lddd, "\spec(\beta_+)"'] \arrow[dd, "\spec \beta" description, dashed] &                       & {\bf Y} \times_{\bf X} \spec A \arrow[ld, "\pi_1"] \arrow[dd, "\pi_2"] \arrow[lddd, "\spec(\iota) \circ \pi_2" description] \\
		&                                                                                                                                   & \bf Y \arrow[d, "f"'] &                                                                                                                 \\
		& \spec A \arrow[r, "i"'] \arrow[ld,"\spec(\iota)"]                                                                                                & \bf X                 & \spec A, \arrow[l, "i"] \arrow[ld, "\spec(\iota)"]                                                               \\
		\spec A_+ &                                                                                                                                   & \spec A_+ \arrow[ll,equal]  &                                                                                                                
	\end{tikzcd}
\end{equation}where $j$ is a Zariski open inclusion (\'{e}tale, resp.) if $i$ is a Zariski open inclusion (\'{e}tale, resp.). 
Thus, we  get the sub-diagram
\begin{equation} 
\begin{tikzcd}
	\spec B \arrow[d, "\spec \beta"'] \arrow[rr, "j:=\pi_1 \circ j'", hook] &  & \bf Y \arrow[d, "f"] \\
	\spec A \arrow[rr, "i", hook]                                           &  & \bf X,               
\end{tikzcd}
\end{equation} which provides the desired local representative, $\spec \beta$, for the morphism $f: \bf Y \rightarrow X$.
\paragraph{(2) Isotropic structures $\varLambda' \text{ and } d_{dR}\varLambda'.$} Since locally $f^*\alpha_0 \sim 0$ via $ \varLambda' $, there is a corresponding element, say $-\varLambda \in (\Omega_B^1)^{k-1}$, with $d_B(-\varLambda) = \beta_* (\alpha_0)$. Letting $\omega_0:=(d_{dR}\alpha_0, 0, \dots) \in \mathcal{A}^{(2,cl)}(\spec A, k)$,  $ d_{dR}\varLambda' $ lifts to an isotropic structure $ d_{dR}\varLambda $ on $\spec \beta :\spec B \rightarrow (\spec A, \omega_0)$, because one gets $d_B(d_{dR}\varLambda)= \beta_* (d_{dR}\alpha_0) \text{ and } d_{dR}(d_{dR}\varLambda)=0$. Moreover, we simply get $\rho:= \rho'|_{\mathbb{T}_B} : \mathbb{T}_B \rightarrow f^*(\ker \alpha_0)$ commuting with all the structure as $j: \spec B \hookrightarrow \bf Y$ is a Zariski open inclusion. %, where $d_{dR}\varLambda=(d_{dR}\varLambda, 0, \dots)$. %with $d_{dR}\varLambda^i \in (\Lambda^{2+i} \Omega_B^1)^{k-1+i}$, such that $d(d_{dR}\varLambda^0) = \beta_* (d_{dR}\alpha_0)$ and $d(d_{dR}\varLambda^{i+1})=0$ for $i\geq 0.$ 

Applying Lemma \ref{lemma_Lag simplification} to the isotropic structure $ d_{dR}\varLambda $, find $G'\in B^k$ and $\psi\in (\Omega_B^1)^{k-1}$ satisfying  \begin{equation} \label{defining eqs with G prime}
	dG'=-\beta(H+H_+), \ d_{dR}G'+d\psi=-\beta_*(\phi+\phi_+)
\end{equation}  such that the isotropic structure $d_{dR}\varLambda$ is homotopic to $  (\frac{1}{k-1}d_{dR}\psi, 0, 0, \dots) $. For the proof, we refer to \cite[Prop. 4.1]{JS}. Comparing with (\ref{defining relations for B,d_B}), we simply have $G', \frac{1}{k-1}\psi$ in place 
of $G, \psi.$%Hence, we can write $(k-1)\varLambda^0=\psi + d_{dR}\theta$ for $\theta \in B^{k-1}$. For the rest of the proof, we will use this simplification. %That is, we consider $ (k-1)d_{dR}\varLambda^0 = d_{dR}\psi$ satisfying $ dG'=-\beta(H+H_+) \text{ and } d_{dR}G'+d\psi=-\beta_*(\phi+\phi_+) $ such that $ \varLambda^0 $ is the non-degenerate isotropic structure on the morphism $\spec \beta :\spec B \rightarrow (\spec A, \alpha_0)$.
\paragraph{(3) Choosing coordinates.} %In addition to the initial assumption $k\not\equiv 3\mod 4$, 
For simplicity, assume further that $k$ is odd with $k\equiv 1 \mod 4$ so that $\ell=-\lfloor(k+1)/2\rfloor$ and $s:=-\lfloor k/2\rfloor$. So far, we obtained a submersion $\beta_+: A_+ \rightarrow B$ minimal at $q$, a smooth $\K$-algebra $A(0)$ of dimension $m_0$ with degree 0 generators $x_1^0, x_2^0, \dots, x_{m_0}^0$, and a set of variables $x_j^{-i}, y_j^{k+i},z^k \text{ in } A$  such that $A_+$ is freely generated over $A(0)$ by the variables $x_j^{-i},$ only. Following \cite[Sections 4.2-4.6]{JS}, we have:

\begin{itemize}
	\item Since $\beta_+^0: A^0 \rightarrow B^0$ is smooth, localizing $B^0$ if necessary, there exist  $u^0_1, \dots, u^0_{n_0}$ in $B^0$ such that $d_{dR}\tx_{1}^0, \dots, d_{dR}\tx^0_{m_0}, d_{dR}u^0_1, \dots, d_{dR}u^0_{n_0}$ form a $ {B^0}$-basis for $ \Omega_{B^0}^1 $   with $\tx^0_j=\beta^0_+(x^0_j)$.
	\item Since $\beta_+: A_+\rightarrow B$ is a submersion of standard form cdgas (and  $A_+$ is freely generated over $A^0$ by the variables $x_j^{-i},$ only),  $B$ is a free graded algebra over $B(0)$ generated by the variables  $\{\tx^{-i}                                                                                                                 ç_j=\beta_+(x^{-i}_j)\}$ and some other set of generators in degrees $-1,-2,\dots, (k-1)$ as $\Omega_B^1$ has Tor-amplitude in $[k-1,0].$ 
\end{itemize}
Therefore, we can consider the cdga $B$, as a \emph{commutative graded algebra}, to be the free graded algebra over $B(0)$ generated by the variables  
\begin{align} \label{generators of B}
	& \tx_1^{-i}, \tx_2^{-i}, \dots, \tx_{m_i}^{-i}  \text{ in degree } (-i) \text{ for } i= 1, 2, \dots, \ell, \nonumber \\
	&u_1^{-i}, u_2^{-i}, \dots, u_{n_i}^{-i} \text{ in degree } (-i) \text{ for } i=1,\cdots s, \nonumber \\
	&v_1^{k-1+i}, v_2^{k-1+i}, \dots, v_{n'_i}^{k-1+i} \text{ in degree } (k-1+i), \text{ for } i=0,1,\cdots s,
\end{align} where $\tx^{-i}_j=\beta_+(x^{-i}_j)$. The full descriptions of $\beta$ and  $d_B$ will be discuss below.

Let us first verify that we must have $n_i=n'_i$ in (\ref{generators of B}). Recall that the  complex 
\begin{equation}
	\Bigg( \big(\mathbb{T}_{B} \otimes_{B}H^0(B)\big) \oplus   \big(\ker \alpha_0 [-1]\otimes_{A}H^0(B)\big) ,       \begin{pmatrix}
		d_B & 0 \\
		\rho & d
	\end{pmatrix}   \Bigg)
\end{equation}serves as a model for $ \mathbb{T}_{\rho} \otimes_{B}H^0(B)$. So, we can write $ \mathbb{T}_{\rho} \otimes_{B}H^0(B)$ and $ \Omega_B^1[k-1] \otimes_{B}H^0(B) $ as direct sums of $H^0(B)$-modules. For each $i\geq 0$, we have:
\begin{align*}
	\left(\mathbb{T}_{\rho} \otimes_{B}H^0(B)\right)^i&= \langle \partial/\partial \tx^{-i}_j\rangle \oplus \langle \partial/\partial u_j^{-i}, \partial/\partial v_{j}^{-i} \rangle \oplus
	\ \langle \partial/\partial x^{-i+1}_j\rangle \oplus \langle \partial/\partial y^{1-i}_j\rangle, \nonumber \\
	\left(\Omega_B^1[k-1] \otimes_{B}H^0(B)\right)^{-i}&= \langle d_{dR}\tx^{k-i+1}_j\rangle\oplus \langle d_{dR}u_j^{k-1+i}, d_{dR} v_j^{k-1+i}\rangle.
\end{align*} By the non-degeneracy condition, $\chi_{\rho}: 	\mathbb{T}_{\rho} \rightarrow \Omega_B^1[k-1]$ is a quasi-isomorphism, where the map $ \chi_{\rho}\otimes_{B}H^0(B) $ is given degree-wise by 
\begin{equation}
\chi^{-*}_{\rho}\otimes_{B}H^0(B)=\begin{pmatrix}
	d_{dR}\varLambda \ \cdot & d_{dR}\varLambda \ \cdot & \star & \beta_* \circ d_{dR}\alpha_0|_{\ker \alpha_0} \ \cdot \\
	d_{dR}\varLambda \ \cdot &  d_{dR}\varLambda \ \cdot & \star & 0
\end{pmatrix}.
\end{equation} Here, $ \star$'s denote some arbitrary morphisms, $\beta_* \circ d_{dR}\alpha_0 \ \cdot: \partial/\partial y^{\mu}_j \mapsto d_{dR}x^{k-\mu}_j \mapsto d_{dR}\tx^{k-\mu}_j$, and $d_{dR}\varLambda \ \cdot$ maps $\langle \partial/\partial u_j^{-i}, \partial/\partial v_j^{-i} \rangle \mapsto \langle d_{dR}u_j^{k-1+i},d_{dR} v_j^{k-1+i} \rangle$, all up to sign. 
Also due to the non-degeneracy of contact isotropic structure, the map $ f^*(d_{dR}\alpha_0|_{\ker \alpha_0}) \ \cdot : f^*(\ker \alpha_0)\rightarrow f^*(\ker^{\vee} \alpha_0[k])$ is a quasi-isomorphism as well, see Observation \ref{observation_non-deg isot str induces weak eqv on pullback}.

From the proof of \cite[Theorem 3.7]{JS}, this non-degeneracy condition holds true if and only if  the morphism $\chi^{-i}_{\rho}\otimes_{B}H^0(B): \langle \partial/\partial u_j^{-i}, \partial/\partial v_j^{-i} \rangle \mapsto \langle d_{dR}u_j^{k-1+i},d_{dR} v_j^{k-1+i} \rangle$ is an isomorphism for each degree. In that case, the map also reduces to
\begin{equation}
	\chi_{\rho}^{*}=\begin{pmatrix}
		0 & 0 & \star & \beta_* \circ d_{dR}\alpha_0 |_{\ker \alpha_0}\ \cdot \\
		0 &  d_{dR}\varLambda \ \cdot & \star & 0
	\end{pmatrix}.
\end{equation} %\footnote{It means $ (\chi^{-i}|_{\ker \varLambda^0})|_q $ is an isomorphism of $\mathbb{K}$-vector spaces for $q\in \spec H^0(B)$ and for each $i$.}.  That is, $\chi$ maps $ \ker \varLambda^0 $ isomorphically onto its image. 
Therefore,  we conclude that $n_j=n'_j$ in (\ref{generators of B}), % due to the non-degeneracy condition, and we 
which also implies (using the isomorphisms $ \chi^{-i}_{\rho}\otimes_{B}H^0(B) $ above) that we have the form $$d_{dR}\varLambda= \sum_{i=0}^{s} \sum_{j=1}^{n_i} d_{dR}u_j^{-i}d_{dR}v_j^{k-1+i}  \quad \text{such that} \quad  (k-1)d_{dR}\varLambda= d_{dR}\psi.$$
Additionally, it follows from the observations above that the map $\rho:= \rho'|_{\mathbb{T}_B} : \mathbb{T}_B \rightarrow f^*(\ker \alpha_0)$, when restricted to $q\in \spec H^0(B)$, is given as $\rho\otimes_{B}H^0(B):\mathbb{T}_B\otimes_{B}H^0(B)\rightarrow \ker \alpha_0\otimes_{A}H^0(B),$  sending $\langle \partial/\partial \tx^{-i}_j \rangle_{H^0(B)} \mapsto \langle \partial/\partial x^{-i}_j \rangle_{H^0(B)}$, with trivial action on the other vector fields. Therefore, the pair $(\varLambda, \rho)$ is identified with the one in Construction \ref{Leg model}.

Now,  our goal is to find an element $ \varLambda \in \Omega_B^1[k-1]$ satisfying the desired properties above. Define such $\varLambda $ to be the element of the form 
\begin{equation}
	\varLambda = \sum_{i=0}^{s} \sum_{j=1}^{n_i}  u_j^{-i}\dR v_j^{k-1+i}.
\end{equation}
Joyce and Safronov showed in  (their proof of) \cite[Theorem 3.4]{JS} that  one can write $\psi$ explicitly\footnote{They arrived at the desired formulas after using suitable homotopic replacement for the morphism $\beta$, but keeping the remaining data fixed. For details, we refer to \cite[Section 4.5.]{JS}} as%\footnote{As before, we may use a homotopic replacement, if needed.} 
\begin{equation}
\psi:= \displaystyle \sum_{i=0}^{s} \sum_{j=1}^{n_i} \big[-iu_j^{-i}d_{dR}v_j^{k-1+i} +(-1)^{(1-i)k}(k-1+i)v_j^{k+i}d_{dR}u_j^{-i}\big].
\end{equation} Therefore, we observe
\begin{equation} \label{11111}
	(k-1)\varLambda=\psi + \dR \Big[\sum_{i,j} (-1)^{(1-k-i)} (1-k-i) (-1)^{(1-i)k}v_j^{k+i-1}u_j^{-i}\Big],
\end{equation}which will be useful to obtain the desired representations as in Construction \ref{Leg model}.  

\paragraph{(4) The morphism $\beta: A\rightarrow B$.} Expanding $ d_{dR}G'+d\psi=-\beta_*(\phi+\phi_+) $ with the  representations above and using $\tx^{-i}_j=\beta_+(x^{-i}_j)$, Joyce and Safronov show  \cite[Theorem 3.7]{JS} that the action of $\beta$ is given by \[\beta|_{A_+}=\beta_+ \quad \text{and} \quad \beta(y_j^{k+i})= (-1)^{-i+1} \partial G/\partial \tx_j^{-i}. \]
On the other hand, using (\ref{defining eqs with G prime}), (\ref{11111}), and $d_B(-\varLambda) = \beta_* (\alpha_0)$, we  compute  
\begin{align*}
	(1-k)\beta_*(\alpha_0)&=(k-1) d_B\varLambda \\
	&= d_B\psi + d_B\circ \dR [\cdots] \\
	&= -\dR G' - \beta_*(\phi + \phi_+) + d_B\circ \dR [\cdots]\\
	&= -\dR G' +(1-k) \beta_*(\alpha_0) \ + (k-1) \dR \beta (z^k) \\ &- \dR\Big[\sum_{i,j}  i \tx_j^{-i}\partial G'/\partial \tx_j^{-i} \Big] + {d_B\circ \dR [\cdots]},
\end{align*} which implies that the action of $\beta$ on $  z^k $ is determined by \begin{align*}
(k-1)\beta(z^k)&:=G'+\sum_{i,j} i \tx^{-i}_j \partial G'/\partial \tx_j^{-i} \\ &+ d_B\Big[\sum_{i,j} (-1)^{(1-k-i)} (1-k-i) (-1)^{(1-i)k}v_j^{k+i-1}u_j^{-i}\Big].
\end{align*} Thus, the action of $\beta$ on $x^{-i}_j, y_j^{k+i}, z^k \in A$ is given as  in Construction \ref{Leg model}.

\paragraph{(5) The differential $d_B$ and classical master equation.} Expanding $ d_{dR}G'+d\psi=-\beta_*(\phi+\phi_+) $ with the  representations above and using $\tx^{-i}_j=\beta_+(x^{-i}_j)$, \cite[Theorem 3.7]{JS} also proves that 
the differential on $B$ is can be written as in (\ref{differantial on B}). Namely, we get  the formulas \begin{align*}
	&d_B|_{B^0}=0, \quad d_B\tx_j^{-i}=(-1)^{1-i} \beta_+(H_j^{-i+1}), \quad \text{ for } i= 1, 2, \dots, \ell, \\ &d_Bu_j^{-i}=(-1)^{(1-i)k} \partial G'/ \partial v_j^{k-1+i}, \quad d_Bv_j^{k-1+i}= \partial G'/ \partial u_j^{-i} \quad \text{ for } i=0, 1, 2, \dots, s,
\end{align*} which correspond to the ones in Construction \ref{Leg model}.

Likewise, using the formula of the differential above to expand $ d_BG'=-\beta(H+H_+) $ also implies that $G'$ must satisfy the \textit{relative version} of  classical master equation (\ref{relative CME}).

We then conclude that the data of $(A, \alpha_0, B,d_B,\beta, \varLambda, \rho)$ is identified with the one in Construction \ref{Leg model}, which (almost) completes the proof.

\paragraph{Other cases of $k.$} As mentioned before, the other cases of $k$ can be studied using similar arguments with some modifications. The first part of the proof above will remain unchanged. When it comes to choosing coordinates and writing down $\psi$ explicitly, different sets of variables must be used. For details, we refer to \cite[Section 4]{JS}. It is straightforward to apply these modifications to our case, although it can sometimes be cumbersome. We leave the details to the reader. 

We then complete the proof of Theorem \ref{Leg model}, and hence that of Theorem \ref{THM3}.

\epf

\section{Concluding remarks}
In this work,  we introduce \emph{Legendrians}, or more precisely, \emph{Legendrian structures} in the derived context. Formal definitions along with some relevant results are provided for \emph{derived (Artin) stacks}, see Definitions \ref{defn_isotropic for Leg} \& \ref{defn_LEG strc},   Proposition \ref{prop_Leg on f with point target}, and Theorem \ref{thm_on zero section}.

On the other hand, Construction \ref{Leg model} describes the prototype affine model for Legendrian structures on the morphisms of \emph{derived affine schemes}; Theorem \ref{Legendrian-Darboux} provides a Legendrian neighborhood theorem for the Legendrians in \emph{contact derived schemes}. Following the similar modifications in \cite{Berktav2024}, we may argue that Theorem \ref{Legendrian-Darboux} might still be valid for Legendrians in \emph{contact derived Artin stacks}. However, this statement needs further clarification and formal treatment, which we do not intend to provide. We leave this task to curious readers.
%----------------------------------------------------
 \newpage
 
 \appendix
 
 \section{Some  derived algebraic geometry} \label{section_prelim}
 In this section, we  review  derived algebraic geometry (DAG) and mention some useful results. We primarily follow \cite{Brav,BenBassat,JS}. For more technical details on DAG, we refer to \cite{Toen,ToenHAG,Luriethesis}. 
 
 %In this paper, we essentially use the functorial approach to define (higher) spaces of interest. 
 \paragraph{Higher spaces as homotopy sheaves.} Recall that using Yoneda's embedding,  spaces can be thought of as \textit{certain functors} in addition to the standard ringed-space formulation, as outlined in the diagram  \cite{Vezz2} 
 \begin{equation} \label{the diagram}
 	\begin{tikzpicture}
 		\matrix (m) [matrix of math nodes,row sep=1.5em,column sep=6.5em,minimum width=1.5 em] {
 			CAlg_{\mathbb{K}}   & Sets  \\
 			&  Grpds \\
 			cdga_{\mathbb{K}}^{\leq 0} &  Ssets. \\};
 		\path[-stealth]
 		(m-1-1) edge  node [left] { } (m-3-1)
 		edge  node [above] {{\small schemes}} (m-1-2)
 		(m-1-1) edge  node [below] {} node [below] {{\small stacks}} (m-2-2)
 		(m-1-1) edge  node [below] {} node [below] {{\small  higher stacks}} (m-3-2)
 		(m-3-1) edge  node  [below] {{\small derived stacks}} (m-3-2)
 		
 		(m-1-2) edge  node [right] { } (m-2-2)
 		(m-2-2) edge  node [right] { } (m-3-2);
 		%edge [dashed,-] (m-2-1);
 	\end{tikzpicture}
 \end{equation} Here $ CAlg_{\mathbb{K}} $ denotes the \emph{category of commutative $\K$-algebras}, and $ cdga_{\K}^{\leq 0}$ is the category\footnote{We actually mean the \emph{$\infty$-category} associated to the model category $ cdga_{\K} $,  with its natural model structure
 for which equivalences are quasi-isomorphisms, and fibrations are
 epimorphisms in strictly negative degrees.} of \emph{commutative differential graded $\mathbb{K}$-algebras in non-positive degrees}. Denote by $St_{\K}$ the $\infty$-category of (higher) $\K$-stacks, where objects in $St_{\K}$ can be ``defined" via  Diagram \ref{the diagram} above as  \emph{homotopy sheaves}. 
 
 Recall that, in classical algebraic geometry, we  have the  ``spectrum functor" 
 \begin{equation*}
 	\spec : (CAlg_{\mathbb{K}})^{op} \rightarrow St_{\K}.
 \end{equation*} We then call an object $X$ of $ St_{\K} $ an \bfem{affine $\K$-scheme} if $X\simeq \spec A$ for some $A\in CAlg_{\mathbb{K}}$; and a \bfem{$\K$-scheme} if it has an open cover by affine $\K$-schemes. Here, $\spec A$ is the usual \emph{ prime spectrum} of $A\in CAlg_{\mathbb{K}}$. 

In derived geometry, there also exists an appropriate concept of a \bfem{spectrum functor} described as the \emph{right adjoint to the global algebra of functions functor} 
\begin{equation*}
	\Gamma: dSt_{\K} \leftrightarrows (cdga_{\K}^{\leq 0})^{\mathrm{op}}: \spec,
\end{equation*} where $dSt_{\K}$ denotes the \bfem{$\infty$-category of derived stacks}, with objects being \emph{$\infty$-sheaves on the site $(dAff)^{\mathrm{op}}:= cdga_{\mathbb{K}}^{\leq 0}$}. We formally denote elements in this opposite category $ (cdga_{\mathbb{K}}^{\leq 0})^{\mathrm{op}} $ by $\spec A$. 

\begin{definition}
	An object ${X}$ in $dSt_{\mathbb{K}}$ is called an \bfem{affine derived $\mathbb{K}$-scheme} if  $X\simeq \spec A $ for some cdga $A \in cdga_{\K}^{\leq 0}$. An object $ X$ in $dSt_{\mathbb{K}}$ is then called a \bfem{derived $\mathbb{K}$-scheme} if it can be covered by Zariski open affine derived $\mathbb{K}$-schemes.
\end{definition}
Denote by $dSch_{\K} \subset dSt_{\mathbb{K}}$ the full \bfem{$\infty$-subcategory of  derived $\mathbb{K}$-schemes}, and we simply write $dAff_{\K} \subset dSch_{\mathbb{K}}$ for the full \bfem{$\infty$-subcategory of  affine derived $\mathbb{K}$-schemes.} Note that $\spec: (cdga_{\K}^{\leq 0})^{\mathrm{op}} \rightarrow dAff_{\K} $ gives an equivalence of $\infty$-categories.

\begin{remark}
	\begin{enumerate}
		\itemsep=0.2cm
		\item [ ]
		\item 	The notation $\spec (-)$ for cdgas must be understood in a  \emph{categorical sense}, meaning we do not use any explicit constructions like prime spectrum of a ring (or locally ringed space) as in ordinary algebraic geometry. 
		\item The choice of $\spec$-notation remains meaningful in geometric terms because the points of a derived scheme correspond with those of its truncation, which is simply an ordinary scheme. That is, for a cdga $A$, the points of $\spec A$ are just the ones in $\spec (H^0(A))$, which is the ordinary prime spectrum of $H^0(A)\in CAlg_{\mathbb{K}}$. 
	\end{enumerate}
\end{remark}

More generally, we  have the following definitions, extending derived schemes.

\begin{definition}
	A \bfem{derived stack} $ {X} $ is  a  functor  $$ {X}:  cdga_{\K}^{\leq 0} \rightarrow  Grp_{\infty}, \ \ A\mapsto X(A) \simeq Map_{dStk_{\mathbb{K}}}(\spec A, {X}),$$ satisfying a descent condition.
	%That is, objects in $dSt_{\mathbb{K}}$ are \emph{$\infty$-presheaves on the site $(dAff)^{\mathrm{op}}\simeq cdga_{\mathbb{K}}^{\leq 0}$  satisfying a descent condition}. 
	For more details, we refer to \cite{ToenHAG}. 
\end{definition}

\begin{definition}(Formal)
	An object $X \in dSt_{\K}$ is called a \bfem{derived Artin $\K$-stack} if it is $m$-geometric for some $m$, and the underlying classical stack is 1-truncated. For  details, we refer to \cite[$\S$ 5.1]{Luriethesis} or \cite[$\S$ 1.3.3]{ToenHAG}.% More formally, we have
\end{definition}
The upshot is that any such object $ X$ of $ dSt_{\K}$ comes with a smooth surjective morphism  $\varphi: U \rightarrow  X$ with $U$ a  derived $\K$-scheme. We call such morphism an \bfem{atlas.} Therefore, the following definition will be sufficient for our purposes.

%\begin{definition} %(\cite[Def. 1.3.3.1]{ToenHAG} or \cite[Section 5.1]{Luriethesis})
%	An object $X \in dSt_{\K}$ is called a \textit{derived Artin $\K$-stack} if it is $m$-geometric for some $m$.
%\end{definition} 

\begin{definition} (Informal)
	By a \bfem{derived Artin $\K$-stack},  we mean an object $X$ of $ dSt_{\K}$ possessing an atlas (smooth of some relative dimension) near each point of $X$.
\end{definition} Notice that when such atlases are simply Zariski open inclusions (rel. dim 0), one recovers the usual structure of a derived scheme. 

\begin{definition}We   define the \bfem{algebra of functions} on $X \in dSt_{\K}$ by \begin{equation}\label{defn_O(X) as inf limit}
		\mathcal{O}(X)= \displaystyle \lim_{A \in cdga_{\K}^{\leq 0}, \ \spec A \rightarrow X} A.
	\end{equation} Likewise,  the \bfem{$\infty$-category of quasi-coherent sheaves on $X$} is given by 
	\begin{equation}\label{defn_QCoh(X) as inf limit}
		\mathrm{QCoh}(X)= \displaystyle \lim_{A \in cdga_{\K}^{\leq 0}, \ \spec A \rightarrow X} \mathrm{Mod}_A,
	\end{equation}where we take the limit in the $\infty$-category of cdgas.
\end{definition}

 %We should again note that throughout this paper, $ \mathbb{K} $ will be an algebraically closed field of characteristic zero. We also assume that all classical $ \K $-schemes are \emph{locally of finite type}, and all derived $ \K $-schemes $ \textbf{X} $ are  \emph{locally finitely presented}, by which we mean that $ \textbf{X} $ can be covered by Zariski open affines $\spec A,$ where $A$ is a finetely presented cdga over $\K.$

 %In addition to the spectrum functors $\spec, \cspec$ above, there is a natural \emph{truncation functor} $\tau: dSt_{\K} \rightarrow St_{\K}$, along with a fully faithfull left adjoint \emph{inclusion functor} $\iota: St_{\K} \hookrightarrow dSt_{\K}$, which can be thought of as an embedding of classical algebraic $\K$-spaces into derived spaces. In this regard, for a cdga $A$ there exists an equivalence $\tau \circ \spec A \simeq \cspec H^0(A)$. This means that if $\bf X$ is a (affine) derived $\K$-scheme, then its truncation $X=\tau(\textbf{X})$ is a (affine)  $\K$-scheme. Therefore, we can consider a derived $\K$-scheme $\bf X$ as \emph{an infinitesimal thickening of its truncation} $X$. It  follows that points of a derived $\K$-scheme $\bf X$ are the same as points of  of its truncation $X$. So the main difference between $X$ and $\bf X$ is in fact encoded by the scheme structure, not by the points!
 
 \paragraph{Nice local models for derived $\K$-schemes.}  Let us first recall some basic concepts from \cite{Brav}.
 
 \begin{definition} \label{defn_standard forms}
 	$ A \in cdga_{\mathbb{K}}^{\leq 0}$ is of \bfem{standard form} if  $A^0$ is a smooth finitely generated $\mathbb{K}$-algebra; the  module $\Omega^1_{A^0}$ of K\"{a}hler differentials is free $A^0$-module of finite rank; and the graded algebra $A$ is freely generated over $A^0$ by finitely many generators, all in negative degrees.
 	
 \end{definition}    
 
 In fact, there is a systematic way of constructing such cdgas starting from a smooth $\K$-algebra $A^0:=A(0)$ via the use of a sequence of localizations. More precisely, for any given $n\in \N$, we can inductively construct a sequence of cdgas 
 \begin{equation} \label{A(n) construction}
 	A(0) \rightarrow A(1) \rightarrow \cdots \rightarrow A(i)\rightarrow \cdots A(n)=:A,
 \end{equation} where  $ A^0:=A(0) $, and $A(i)$ is obtained from $A(i-1)$ by adjoining generators in degree $-i$, given by $M^{-i}$, for all $i$. Here, each $M^{-i}$ is a free finite rank module (of degree $-i$ generators) over $A(i-1)$. Therefore, the underlying commutative graded algebra of $A=A(n)$ is freely generated over $A(0)$ by finitely many generators, all in negative degrees $-1,-2,\dots, -n$. For more details, we refer to \cite[Example 2.8]{Brav}.

 \begin{definition} \label{1st defn of minimalty}
 	A standard form cdga $ A $ is said to be \bfem{minimal} at  $p \in \cspec H^0(A)$ if $ A=A(n) $ is defined by using the minimal possible  numbers of graded generators in each degree $\leq 0$ compared to all other cdgas locally equivalent to  $ A $ near $p$. %(There will be an equivalent definition below, see Definition \ref{2nd defn of minimalty}.)
 \end{definition}  
 
 \begin{definition}
 	Let $A$ be a standard form cdga.  $A'\in cdga_{\mathbb{K}}^{\leq 0}$ is called a \bfem{localization} of $A$ if $A'$ is obtained from $A$ by inverting an element $f\in A^0$, by which we mean $A'=A\otimes_{A^0}A^0[f^{-1}]$. %for $f\in A^0$. 
 	
 	$A'$ is then of standard form with $A'^0 \simeq A^0[f]$. If $p\in \cspec H^0(A)$ with $ f(p)\neq0 $, we say $A'$ is a \emph{localization of $A$ at $p$}.
 \end{definition}
 With these definitions in hand, one has the following observations: 
 \begin{observation}
 	Let $A$ be a standard form cdga. If $A'$ is  a \emph{localization} of $A$, then $\spec A' \subset \spec A$ is a Zariski open subset. Likewise, if $A'$ is a \emph{localization of $A$ at $p\in \spec H^0(A)\simeq \tau(\spec A)$}, then  $\spec A' \subset \spec A$ is a Zariski open neighborhood of $p.$
 \end{observation}
 \begin{observation}
 	Let  $ A$ be \emph{a standard form} cdga, then there exist generators $x_1^{-i}, x_2^{-i}, \cdots, x_{m_i}^{-i} $ in $A^{-i}$ (after localization, if necessary) with $i= 1, 2, \cdots, k$  \ and $m_i \in \mathbb{Z}_{\geq 0}$ such that 
 	\begin{equation}
 		A = A(0) \big[x_j^{-i} : i= 1, 2, \dots, k, \ j= 1,2, \dots, m_i\big],
 	\end{equation} where the subscript $j$ in $x_j^i$ labels the generators, and the superscript $i$ indicates the degree of the corresponding element.  So, we can consider  $A$ as a \emph{graded polynomial algebra over $A(0)$ on finitely many generators, all in negative degrees.} 
 \end{observation}%Moreover, since $A^0$ is finitely generated $\mathbb{K}$-algebra, there exists a surjection $\mathbb{K}[x_1^0, \cdots, x_{m_0}^0] \twoheadrightarrow A^0$ such that $A^0= \mathbb{K}[x_1^0,x_2^0, \cdots, x_{m_0}^0]/ I$ for some ideal $I$. %We shall be interested in the cases where the ideals $I$ are also finitely generated. Then $A^0$ is indeed said to be \emph{finitely presented}.
 %	\item In general, we are interested in minimal standard form cdgas for which each number $m_i$ above is taken to be the least possible number such that $m_i=dim (H^i (\mathbb{L}_{A} |_p))$  and $d^i |_p =0$ for $i= -1, -2, \cdots, k$, and $p \in spec (H^0(A))$. 
 
 %\begin{definition}
 %	We then define the \emph{virtual dimension} of $A$ to be the integer $\vdim A= \sum _i (-1)^im_i$. 
 %\end{definition}
 
 The following theorem outlines the central result from \cite{Brav} concerning the construction of useful local algebraic models for  derived $\K$-schemes.  The upshot is that given a derived $\K$-scheme $\bf X$ (locally of finite presentation) and a point $x\in \bf X$, one can always find a ``refined"  affine neighborhood $ \spec A $ of $x$, which is very useful for explicit presentations.%; and any two refined neighborhoods of $x$ can be comparable on the overlap.    
 In short, we have:
 \begin{theorem} \cite[Theorem 4.1]{Brav}\label{localmodelthm}
 	Every derived $\mathbb{K}$-scheme $X$ is Zariski locally modelled on $ \spec A $ for a minimal standard form cdga $ A$. %Moreover, if there are two such models $ \spec A \text{ and } \spec B$,	then there exists another minimal standard form cdga $ C $ with a morphim $\spec C \rightarrow {\bf X}$, which relates $ \spec A \text{ to } \spec B$ in certain way.
 \end{theorem}

 \paragraph{Nice local models for cotangent complexes of derived schemes.} Given	$ A \in cdga_{\mathbb{K}}^{\leq 0}$, $ d$ on $A$ induces a differential on $\Omega_A^1$, denoted again by $ d$. This makes $\Omega_A^1$ into a $\rm dg$-module $(\Omega_A^1, d)$ with the property that $\delta \circ  d = d \circ \delta,$ where $\delta: A \rightarrow \Omega_A^1$ is the universal derivation of degree zero. Write  the decomposition of $ \Omega^1_{A} $ into graded pieces 
 $ \Omega^1_{A} =  \bigoplus_{k=-\infty}^0 \big(\Omega^1_{A}\big)^k $
 with the  differential $d: \big(\Omega^1_{A}\big)^k \longrightarrow \big(\Omega^1_{A}\big)^{k+1}$. Then we define the \bfem{de Rham algebra of $A$} as a  double complex 
 \begin{equation}
 	DR(A)= Sym_A(\Omega_A^1[1]) \simeq \bigoplus \limits_{p=0}^{\infty} \bigoplus \limits_{k=-\infty}^0 \big(\Lambda^p \Omega^1_{A}\big)^k [p],
 \end{equation} where the gradings $p, k$ are called the \emph{weight} and the \emph{degree}, respectively. Also, there are two derivations (differentials)  on $DR(A)$, namely the \bfem{internal differential d}$ : \big(\Lambda^p \Omega^1_{A}\big)^k [p] \longrightarrow \big(\Lambda^p \Omega^1_{A}\big)^{k+1} [p] $ and the \bfem{de Rham differential} $d_{dR}:\big(\Lambda^p \Omega^1_{A}\big)^k [p] \longrightarrow \big(\Lambda^{p+1} \Omega^1_{A}\big)^k [p+1]$ such that $d_{tot}= d + d_{dR}$ and  
 	$ d^2=d_{dR}^2=0, \text{\ and \ } d\circ d_{dR} =- d_{dR} \circ d. $
 % Here, one also has the natural multiplication on $ DR(A)$:
 %\begin{equation}
 %	\big(\Lambda^p \Omega^1_{A}\big)^k [p] \times \big(\Lambda^q \Omega^1_{A}\big)^{\ell} [q] \longrightarrow \big(\Lambda^{p+q} \Omega^1_{A}\big)^{k+\ell} [p+q].
% \end{equation} 
 
 % \begin{observation}
 	% 	 It should be noted that the constructions of $\Omega_A^1, \ DR(A)$ depend only on the underlying commutative graded algebra of $A$, not on the differential $d$ on $A.$ 
 	%\end{observation}
 	
Note that even if both $\mathbb{L}_{A} \text{ and } \Omega^1_{A}$ are closely related, the identification of $ \mathbb{L}_{A} $ with  $ \Omega^1_{A} $ is not true for an arbitrary   $ A \in cdga_{\mathbb{K}}^{\leq 0}$ \cite{Brav}. But, when $A=A(n)$ is a standard form cdga, we have the following description for the restriction of the cotangent complex $\mathbb{L}_A $ to $\spec H^0(A)$. In this paper, we only give a brief version. More details and the proof can be found in \cite[Prop. 2.12]{Brav}.
 	
 	\begin{proposition}  \label{proposition_L as a complex of H^0 modules}
 		If $A=A(n)$ with $n\in \N$ is a standard form cdga constructed inductively as in (\ref{A(n) construction}), then the restriction of $\mathbb{L}_A $ to $\spec H^0(A)$ is represented by a finite complex of $H^0(A)$-modules.
\end{proposition}
 
 \paragraph{Nice local models for morphisms of derived schemes.} Recall that given a derived $\K$-scheme $\bf X$ (locally of finite presentation) and a point $x\in \bf X$, one can always find a ``refined"  affine neighborhood $ \spec A $ of $x$ that allows us to make more concrete computations over this neighborhood (cf. Theorem \ref{localmodelthm}). Such models are central to construct Darboux-type local forms. Similarly, there are \emph{ nice local representatives} for morphisms $f: \bf Y \rightarrow X$ in $dSch_{\K}$. In this regard, let us outline some key concepts and results from \cite{JS}.
 
 \begin{definition}
 	A morphism $\beta: A \rightarrow B$ of standard form cdgas is called a \bfem{submersion} if the induced morphism $\beta_*: \Omega_A^1 \otimes_{A} B \rightarrow \Omega_B^1$ is injective in each degree.
 \end{definition}
 
 \begin{remark}
 	Submersions are in fact suitable class of morphisms in $ cdga_{\K}^{\leq 0}$ for doing explicit computations concerning relative cotangent complexes. In this regard, if $\beta: A \rightarrow B$ is a submersion of standard form cdgas, then the module of relative K\"{a}hler differentials $\Omega_{B/A}^1$ is a model for the relative cotangent complex, and hence we take $\mathbb{L}_{B/A}=\Omega_{B/A}^1$. 
 \end{remark}

 \begin{theorem} \label{thm_nice models for morphisms} \cite[Theorem 3.2]{JS}
 	Let $f: \bf Y \rightarrow X$ in $dSch_{\K}$, $y\in \bf Y$, $x\in \bf X$ with $f(y)=x$. Let $(A, i: \spec A \hookrightarrow \bf X),$ $p\in \spec H^0(A)$ be a standard form affine neighborhood of $x$ with $i(p)=x$ (as in Theorem \ref{localmodelthm}). Then there exist a standard form cdga $B$, a point $q\in \spec B,$ a submersion $\beta: A \rightarrow B$ minimal at $q$ with $(\spec \beta) (q)=p$, and a Zariski open inclusion $j:\spec B \hookrightarrow \bf Y$, with $j(q)=y$, such that the following diagram (homotopy) commutes:
 	\begin{equation} 
 		\begin{tikzpicture}
 			\matrix (m) [matrix of math nodes,row sep=1.5em,column sep=4.5em,minimum width=1.5 em] {
 				\spec B   & {\bf Y}  \\
 				\spec A &  {\bf X} \\};
 			\path[-stealth]
 			(m-1-1) edge  node [left] {{\small $ \spec \beta $} } (m-2-1)
 			(m-1-1) edge  node [above] {{\small $ j $}} (m-1-2)
 			%(m-1-1) edge  node [below] {} node [below] {{\small stacks}} (m-2-2)
 			%(m-1-1) edge  node [below] {} node [below] {{\small  higher stacks}} (m-3-2)
 			(m-2-1) edge  node  [below] {{\small $ i $}} (m-2-2)
 			
 			(m-1-2) edge  node [right] { {\small $ f $}} (m-2-2);
 			%edge [dashed,-] (m-2-1);
 		\end{tikzpicture}
 	\end{equation}
 \end{theorem}
 
 In brief, any morphism $f: \bf Y \rightarrow X$ in $dSch_{\K}$ can be locally modeled on $\spec \beta$, with $\beta$ a morphism of standard forms cdgas, which can be considered as a ``map between local charts". I.e., $\spec \beta$ plays the role of a \emph{local representative} for $f$ (with suitable affine neighborhoods both on the source and the target).
 
 Let us now present some other useful result (not from \cite{JS}):
 \begin{proposition}\label{prop_pullback of cocone=cocone of the pullback}
 	Let $f: \bf Y \rightarrow \bf X$ be a morphism of derived Artin stacks locally of finite presentation and $\alpha \in \mathcal{A}^p({\bf X},n)$. Then there are equivalence[s]\begin{equation}
 		Cocone(f^*\alpha)\simeq f^* \big(Cocone(\alpha)\big) \text{ and } Cone(f^*\alpha)\simeq f^* \big(Cone(\alpha)\big).
 	\end{equation}
 \end{proposition}
 \pf
 Firstly, there are two homotopy fiber sequences \begin{equation}
 	Cocone(f^*\alpha)\rightarrow f^*\mathbb{T}_{{\bf X}} \xrightarrow{f^*\alpha} \mathcal{O}_{{\bf Y}}[n] \ \text{ and } 
 	f^*Cocone(\alpha)\rightarrow f^*\mathbb{T}_{{\bf X}} \xrightarrow{f^*\alpha} \mathcal{O}_{{\bf Y}}[n],
 \end{equation}where the latter follows from the pullback of $ Cocone(\alpha)\rightarrow \mathbb{T}_{{\bf X}} \xrightarrow{\alpha} \mathcal{O}_{{\bf X}}[n]. $ Then we get the equivalence of exact sequences
 
 \begin{equation}
 	\begin{tikzpicture}
 		\matrix (m) [matrix of math nodes,row sep=3em,column sep=5 em,minimum width=3 em] {
 			Cocone(f^*\alpha)	& f^*\mathbb{T}_{{\bf X}}  & \mathcal{O}_{{\bf Y}}[n]\\
 			f^*Cocone(\alpha)	& f^*\mathbb{T}_{{\bf X}}  & \mathcal{O}_{{\bf Y}}[n],  \\			
 		};
 		\path[-stealth]
 		
 		%horizontal 1st row
 		(m-1-1) edge  node [above] { $ $} (m-1-2)
 		(m-1-2) edge  node [above] { $ $} (m-1-3)
 		%(m-1-4) edge  node [right] { } (m-1-5)
 		%(m-1-5) edge  node [right] { } (m-1-6)
 		%(m-1-6) edge  node [right] { } (m-1-7)
 		%(m-1-7) edge  node [right] { } (m-1-8)
 		%horizontal 2nd row
 		(m-2-1) edge  node [above] { $  $ } (m-2-2)
 		(m-2-2) edge  node [above] { $   $ } (m-2-3)
 		%(m-2-4) edge  node [right] { } (m-2-5)
 		%(m-2-5) edge  node [right] { } (m-2-6)
 		%(m-2-6) edge  node [right] { } (m-2-7)
 		%	(m-2-7) edge  node [right] { } (m-2-8)
 		%vertical
 		(m-1-1) edge  [dashed,->] node [right] {$ \simeq $ } (m-2-1)
 		(m-1-2) edge  node [right] { $ id $ } (m-2-2)
 		(m-1-3) edge   node [right] {$ id $ } (m-2-3)
 		%(m-1-6) edge  node [right] { } (m-2-6)
 		%(m-1-7) edge  node [right] { } (m-2-7);
 		%edge [dashed,-] (m-2-1)
 		;
 	\end{tikzpicture}
 \end{equation}which gives the desired identification. 
 
 Likewise, using the two homotopy cofiber sequences \begin{equation}
 	f^*\mathbb{T}_{{\bf X}} \xrightarrow{f^*\alpha} \mathcal{O}_{{\bf Y}}[n] \rightarrow Cone(f^*\alpha) \ \text{ and } 
 	f^*\mathbb{T}_{{\bf X}} \xrightarrow{f^*\alpha} \mathcal{O}_{{\bf Y}}[n] \rightarrow f^*Cone(\alpha),
 \end{equation} the result follows.
 \epf
 
 \section{Lagrangian-Darboux theorem for symplectic derived schemes}\label{section_models for LAGNBH thm}
 In this section, we outline local models for Lagrangians in derived symplectic geometry and state a Lagrangian neighborhood theorem. In brief, Joyce and Safronov \cite[Examples 3.3 \& 3.5]{JS} provide explicit local models, called \emph{Lagrangian Darboux forms}, for the cases $(1) \ k<0 \text{ with } k\not\equiv 3\mod 4$ and $(2) \ k<0 \text{ with } k\equiv 3\mod 4$. They also showed in \cite[Theorem 3.7]{JS} that any such structure can be locally modeled on one of Lagrangian Darboux forms. In this regard, we overview local models in the following prototype construction/example.
 
 \begin{construction} \label{Lag model}
 	Consider the case $k<0 \text{ with } k\not\equiv 3\mod 4$. Assume $A, \omega$ are in $k$-shifted symplectic Darboux form in the sense of Theorem \ref{Symplectic darboux}. In brief, we assume that $ A $ is a standard form cdga which  is free algebra over a smooth $\K$-algebra $A(0)$ generated by variables $x^{-i}_j, y^{k+i}_j \in A$ so that $\omega:=(\omega^0, 0, 0, \dots)$, with
 	$  \omega^0= \sum_{i,j} d_{dR}x_j^{-i} d_{dR}y_j^{k+i},$ is a $k$-shifted symplectic form. We will not give all the details on each case of $k$ in this paper. Instead, we refer to \cite[Examples 5.8, 5.9, and 5.10]{Brav}.
 	
 	Let $\ell=-\lfloor(k+1)/2\rfloor$\footnote{$ \ell=-k/2 $ if $k$ is even; otherwise $\ell= -(k+1)/2$. Here, $\lfloor - \rfloor$ denotes the \emph{floor function}.}. Write $A^{-i}_+$ for the sub-cdga of $A$ generated by $x^{-i}_j$'s only, for $i=1,\dots, \ell$, $j=0, \dots, m_i.$ Then we consider $A$ as a freely generated algebra on generators $x^{-i}_1, \dots, x^{-i}_{m_i}$ in $A^{-i}_+\subset A$ for $i=1,\dots, \ell$ and $ y_1^{k+i}, \dots, y_{m_i}^{k+i} $ in $A^{k+i}$ for $i=0,\dots, \ell$.
 	
 	Choose an element $H\in A^{k+1}$, called the \emph{Hamiltonian}, satisfying the \textit{classical master equation} \begin{equation} 
 		\displaystyle \sum_{i=1}^{\ell} \sum_{j=1}^{m_i} \dfrac{\partial H}{\partial x_j^{-i}} \dfrac{\partial H}{\partial y_j^{k+i}}=0 \text{ in } A^{k+2}.\end{equation}
 	Then we define the \bfem{internal differential} on $A$ by $d=0$ on $A(0)$,  and  by  \begin{equation} 
 		dx_j^{-i} =  (-1)^{(1-i)(k+1)}\dfrac{\partial H}{\partial y_j^{k+i}} \ \text{ and } \ dy_j^{k+i} =  \dfrac{\partial H}{\partial x_j^{-i}}.
 	\end{equation}
 	Moreover, we can choose an element $\phi \in (\Omega_A^1)^k$ satisfying $dH=0, \ d_{dR}H + d\phi =0, \text{ and } d_{dR}\phi=k\omega^0.$
 	Explicitly, \ we have \begin{equation} \label{defn_phi}
 		\phi:= \displaystyle \sum_{i=0}^{\ell} \sum_{j=1}^{m_i} \big[-ix_j^{-i}d_{dR}y_j^{k+i} +(-1)^{(1-i)(k+1)}(k+i)y_j^{k+i}d_{dR}x_j^{-i}\big].
 	\end{equation}
 	
 	Notice that due to the degree reasons, $H$ must be at most linear in $ y_j^{k+i}. $ Then decomposing $H$ we write 
 	\begin{equation}
 		H=H_+ +  \sum_{i,j} H^{(-i+1)}_jy_j^{k+i}, \text{ with } H_+\in A_+^{k+1}, \ H^{(-i+1)}_j\in A_+^{-i+1}.
 	\end{equation}Likewise, we can define $\phi_+:= -\sum_{i=0}^{\ell} \sum_{j=1}^{m_i} (-1)^{(-i+1)(k+1)} y_j^{k+i}d_{dR}x_j^{-i}$. Then the classical master equation becomes
 	\begin{align}
 		&\sum_{i=1}^{\ell} \sum_{j=1}^{m_i} (-1)^{-i+1} H_j^{(-i+1)}\dfrac{\partial H_+}{\partial x_j^{-i}}=0 \text{ in } A^{k+2}_+ \nonumber \\
 		&\sum_{i=1}^{-i'+1} \sum_{j=1}^{m_i} (-1)^{-i+1} H_j^{(-i+1)}\dfrac{\partial H_{j'}^{(-i'+1)}}{\partial x_j^{-i}}=0 \text{ in } A^{-i'+2}_+
 	\end{align} for $i'=1, \dots, \ell; \ j'=1, \dots, m_{i'}.$ The action of $d$ can be equivalently written as
 	
 	\begin{equation}
 		dx_j^{-i} =  (-1)^{-i+1} H_j^{(-i+1)} \ \text{ and } \ dy_j^{k+i} =  \dfrac{\partial H_+}{\partial x_j^{-i}} + \sum_{i'=-i+1}^{\ell} \sum_{j'=1}^{m_{i'}} y^{k+i'}_{j'}\dfrac{\partial H_{j'}^{-i'+1}}{\partial x_j^{-i}}.
 	\end{equation}Then some (very long) calculations show that $d_+H=0, \ d_{dR}H_+ + d\phi_+ =0, \text{ and } d_{dR}\phi_+=-\omega^0.$ For details, we refer to \cite[Remark 2.15 \& Example 3.3]{JS}.
 	
 	Let $s:=-\lfloor k/2\rfloor$ so that 
 	\[s = 
 	\begin{cases*} \ell, & $ k $ \mbox{ even}, \\
 		\ell+1, & $ k $ \mbox{ odd}.
 	\end{cases*}
 	\]
 	Then \cite[Example 3.3]{JS} chooses a smooth algebra $B^0:=B(0)$ of dimension $m_0+n_0$ and a morphism $\beta^0:A^0\rightarrow B^0$. Localizing $B^0$ if necessary, assume that there exist elements $u^0_1, \dots, u^0_{n_0}$ in $B^0$ such that $\Omega_{B^0}^1= \mathrm{span}_{B^0} \{d_{dR}\tx_{1}^0, \dots, d_{dR}\tx^0_{m_0}, d_{dR}u^0_1, \dots, d_{dR}u^0_{n_0} \}$ with $\tx^0_j=\beta^0(x^0_j)$. Then we define the cdga $B$ to be the free graded algebra over $B(0)$ generated by the variables  
 	\begin{align} 
 		& \tx_1^{-i}, \tx_2^{-i}, \dots, \tx_{m_i}^{-i}  \text{ in degree } (-i) \text{ for } i= 1, 2, \dots, \ell, \nonumber \\
 		&u_1^{-i}, u_2^{-i}, \dots, u_{n_i}^{-i} \text{ in degree } (-i) \text{ for } i=1,\cdots, s, \nonumber \\
 		&v_1^{k-1+i}, v_2^{k-1+i}, \dots, v_{n_i}^{k-1+i} \text{ in degree } (k-1+i), \text{ for } i=0,1,\cdots, s
 	\end{align} Define a morphism $\beta_+:A_+ \rightarrow B$ of cdgas by $\beta_+|_{A^0}=\beta^0$ and $\tx^{-i}_j=\beta_+(x^{-i}_j)$ for all $i,j.$%$ i= 1, 2, \dots, \ell; j=1,\dots, m_i. $
 	
 	Now choose a superpotential $G \in B^k$ satisfying a \textit{relative version} of the classical master equation, which determines the differential $d$ on $B$. In fact, extending $\beta_+$ to a morphism $\beta:A\rightarrow B$ via $\beta|_{A_+}=\beta_+$ and $\beta(y_j^{k+i})= \pm \partial G/\partial \tx_j^{-i}$ would give the relative CME in a compact form
 	\begin{equation} \label{relative CME}
 		\displaystyle \sum_{i=1}^{s} \sum_{j=1}^{n_i} \dfrac{\partial G}{\partial u_j^{-i}} \dfrac{\partial G}{\partial v_j^{k+i}} + \beta(G)=0 .\end{equation} Then the \bfem{differential} $d$ on $B$ is given by \begin{align} \label{differantial on B}
 		d|_{B^0}&=0, \quad d\tx_j^{-i}=(-1)^{1-i} \beta_+(H_j^{-i+1}),\nonumber \\  du_j^{-i}&=(-1)^{(1-i)k} \partial G/ \partial v_j^{k-1+i}, \quad  dv_j^{k-1+i}=\pm \partial G/ \partial u_j^{-i}. 
 	\end{align}
 	Finally, we define an element in $(\Lambda^2\Omega_B^1)^{k-1}$
 	\begin{equation} 
 		h^0:=\displaystyle \sum_{i=0}^{s} \sum_{j=1}^{n_i} d_{dR}u_j^{-i}d_{dR}v_j^{k-1+i}.\end{equation}
 	
 	Joyce and Safronov \cite[Example 3.3]{JS} show that $h:=(h^0, 0, 0, \dots)$ is a non-degenerate isotropic structure on $\spec \beta: \spec B \rightarrow (\spec A, \omega_{can})$, with the equations $dh^0= \beta_* (\omega^0), \  d_{dR}h^0=0,$ where $\omega_{can}=(\omega^0, 0, \dots)$ is the standard $k$-shifted symplectic structure on $\spec A$ as in  Theorem \ref{Symplectic darboux}.
 	
 	Moreover, there exists an element $\psi\in (\Omega_B^1)^{k-1}$ (viewed as the relative version of $\phi$ above) satisfying the equations 
 	\begin{equation}
 		dG=-\beta(H+H_+), \ d_{dR}G+d\psi=-\beta_*(\phi+\phi_+), \text{ and } d_{dR}\psi=(k-1)h^0.
 	\end{equation} In fact, by \cite{JS}, we can write $\psi$ explicitly as 
 	\begin{equation}\label{defn_psi}
 		\psi:= \displaystyle \sum_{i=0}^{s} \sum_{j=1}^{n_i} \big[-iu_j^{-i}d_{dR}v_j^{k-1+i} +(-1)^{(1-i)k}(k-1+i)v_j^{k-1+i}d_{dR}u_j^{-i}\big].
 	\end{equation}
 \end{construction}
 \begin{definition} \label{defn_Lag Darboux forms}
 	Let $ \beta\in \Hom(A,B), A, \omega_{can}, B, \text{ and } h=(h^0, 0, \dots)$ be as above. Then we say that $(A, \omega_{can}, B, \beta, h)$ is of \bfem{Lagrangian Darboux form}.
 \end{definition}
 
 It should be noted that the construction above holds true only for particular values of $k$. The complete treatment requires some modifications, but follows the same logic. For details, see \cite[Examples 3.3, 3.5, 3.6]{JS}.

\bibliographystyle{ieeetr} % Choose Phys. Rev. style for bibliography
%Bibliography styles that can be used instead of prsty are abbrv, alpha, plain, acm, ieeetr, siam and unsrt.
\bibliography{refs}

%\printbibliography %Prints bibliography when using biblatex package

\end{document}